\documentclass[journal,twoside,web]{ieeecolor}
\usepackage{generic}
\usepackage{textcomp}

\usepackage{graphicx}
\usepackage{times}
\usepackage{amsmath,amsbsy,amssymb,mathtools}
\usepackage{color}
\usepackage[us]{datetime}
\usepackage{url}
\usepackage{theorem}
\usepackage{multirow}
\usepackage{footnote}
\usepackage{nomencl}
\usepackage{cite}
\usepackage{algorithm}
\usepackage{algpseudocode}
\usepackage{stackengine}

\usepackage{commath}
\usepackage{dsfont}

\newcommand{\negpar}[1][-1em]{%
  \ifvmode\else\par\fi
  {\parindent=#1\leavevmode}\ignorespaces
}

\newcommand{\R}{\mathbb{R}} 

\theoremstyle{plain}
\newtheorem{theorem}{Theorem}

\newtheorem{lemma}{Lemma}
\newtheorem{corollary}{Corollary}
\newtheorem{remark}{Remark}
\newtheorem{problem}{Problem}
\newtheorem{assumption}{Assumption}
\newtheorem{cvxCondition}{Condition}
\usepackage{lipsum}
\newcommand{\qed}{\hfill\rule{1.5ex}{1.5ex}}

\algnewcommand\algorithmicswitch{\textbf{switch}}
\algnewcommand\algorithmiccase{\textbf{case}}
\algdef{SE}[SWITCH]{Switch}{EndSwitch}[1]{\algorithmicswitch\ #1\ \algorithmicdo}{\algorithmicend\ \algorithmicswitch}%
\algdef{SE}[CASE]{Case}{EndCase}[1]{\algorithmiccase\ #1}{\algorithmicend\ \algorithmiccase}%
\algtext*{EndSwitch}%
\algtext*{EndCase}%

\def\BibTeX{{\rm B\kern-.05em{\sc i\kern-.025em b}\kern-.08em
    T\kern-.1667em\lower.7ex\hbox{E}\kern-.125emX}}

\begin{document}

\onecolumn
\copyright 2021 IEEE. Personal use of this material is permitted. Permission from IEEE must be obtained for all other uses, in any current or future media, including reprinting/republishing this material for advertising or promotional purposes, creating new collective works, for resale or redistribution to servers or lists, or reuse of any copyrighted component of this work in other works.

This work has been submitted to the IEEE for possible publication. Copyright may be transferred without notice, after which this version may no longer be accessible.
\twocolumn
\newpage

\title{
Lax Formulae for Efficiently Solving Two Classes of State-Constrained Optimal Control Problems
}
\author{Donggun Lee, \IEEEmembership{Student Member, IEEE}, and Claire J. Tomlin, \IEEEmembership{Fellow, IEEE} 
\thanks{This research is supported by ONR under the BRC program in multibody control systems, by DARPA under the Assured Autonomy program, and by
NSF grant \#1837244.}
\thanks{Donggun Lee is with the Department of Mechanical Engineering, University of California, Berkeley, USA.
        {\tt\small donggun\_lee@berkeley.edu}}%
\thanks{Claire J. Tomlin is with the Department of Electrical Engineering and Computer Sciences, University of California, Berkeley, USA.
        {\tt\small tomlin@eecs.berkeley.edu}}%
}

\maketitle

\begin{abstract}
This paper presents Lax formulae for solving the following optimal control problems: minimize the maximum (or the minimum) cost over a time horizon, while satisfying a state constraint.
We present a viscosity theory, and by applying the theory to the Hamilton-Jacobi (HJ) equations, these Lax formulae are derived.
A numerical algorithm for the Lax formulae is presented: under certain conditions, this algorithm's computational complexity is polynomial in the dimension of the state.
For each class of optimal control problem, an example demonstrates the use and performance of the Lax formulae.
\end{abstract}

\section{Introduction}

This paper presents Lax formulae for two classes of optimal control problems, in which an optimal control signal is specified to minimize the maximum (or the minimum) cost over a time horizon subject to a state constraint. 
In the two classes of problems, terminal time maximizes (or minimizes) cost, and both a control signal and the terminal time must be determined.
The two classes of problems can be solved by the Hamilton-Jacobi (HJ) equations in \cite{donggunTAC_part1} and \cite{donggun2020CDC}.
The HJ equations are, in general, computed by grid-based methods (such as level-set methods and fast marching method \cite{Osher02}), which require spatial discretization, and computational complexity is exponential in the dimension of the system's state.
It is intractable to compute the solution to the HJ equations using the grid-based methods for high dimensional systems.

Lax and Hopf theory specify an optimization problem that solves HJ equations under particular conditions \cite{hopf1965generalized,bardi1984hopf,chow2019algorithm,donggunTAC_submitted}.
The optimization problem from Lax and Hopf theory is computed by gradient-based methods based on temporal discretization, which allows the computational complexity to be polynomial in the dimension of the state.
Lax formulae assume that the Hamiltonian is convex in the costate, on the other hand, Hopf formulae assume that the terminal cost is convex in the state.
Lax and Hopf formulae have been investigated for optimal control problems or zero-sum games, where the cost is the integration of the stage cost and the terminal cost, and the terminal time is a fixed constant.
On the other hand, in the two classes of problems in this paper, the terminal time is not fixed but need to be specified.

The initial version of the Lax formula and Hopf formula were presented in \cite{lax1957hyperbolic,evans10,bardi1984hopf} for particular state-unconstrained problems: in these, the Hamiltonian only depends on the costate, implying the dynamics of the system and the stage cost only depend on the control input.
Bardi and Evans \cite{bardi1984hopf} first proved these initial version of the Lax and Hopf formulae in a sense of viscosity theory of partial differential equations (PDEs).
The Hopf formula has been generalized to the time-dependent Hamiltonian \cite{lions1986hopf} and linear dynamic systems \cite{darbon2016algorithms}.
Later, Chow et al. \cite{chow2019algorithm} generalized Lax and Hopf formulae under the assumption that the stage cost is convex, and the dynamics is affine in the state and the control. 
Also, the authors' previous paper \cite{donggunTAC_submitted} presents a Lax formula dealing with state constraints for fixed-horizon optimal control problems.

In these prior works, there have been two approaches for deriving the Lax and Hopf formulae.
First, Bardi and Evans \cite{bardi1984hopf} followed two steps: in step 1, an equivalent optimal control problem was found by analyzing the corresponding HJ equation, in which a control signal has to be specified, and in step 2, the step-1 Lax or Hopf formulae is simplified to optimization problems (the step-2 formulae) for which the dimension of the decision variable is the same as that of the state.
For the class of problems Bardi and Evans solved, it turns out that the step-1 Lax and Hopf formulae are convex in the control, even if the given optimal control problem is non-convex in the control, leading to fewer variables to be determined for the step-2 formulae in comparison to the step-1 formulae.
On the other hand, Chow et al. \cite{chow2019algorithm} utilized the strong duality technique to directly find approximate step-2 Lax and Hopf formulae, which shows the connection between HJ PDEs and Pontryagin minimum (or maximum) principle.

This paper proposes a step-1 Lax formula for state-constrained optimal control problems for the two classes of problems.
We leave the proposal of a step-2 Lax formula for future work.
The organization of this paper is as follows. 
Section \ref{sec:probDefinition} defines the two classes of problems. 
Section \ref{sec:HJeq} presents the HJ equations corresponding to the both classes of problems.
Section \ref{sec:Laxformula} proposes a viscosity theory, which is utilized to propose and prove the Lax formulae for the two classes of problems.
Section \ref{sec:cvxAnalysis} presents convexity conditions for the proposed Lax formula under which most gradient-based methods guarantee optimality with polynomial complexity.
Section \ref{sec:NumericalAlg} presents a numerical algorithm to compute a solution to the given two classes of problems by utilizing the proposed Lax formulae.
Section \ref{sec:example} presents two examples to demonstrate the utility and performance of the proposed Lax formulae for each class of problems.
Section \ref{sec:conclusion} concludes this paper.

\section{State-Constrained Optimal Control Problem and the Hamilton-Jacobi Equation}
\label{sec:probDefinition}

In this paper, we consider the state trajectory ($\mathrm{x}$) solving the following ordinary differential equation (ODE):
\begin{align}
    \dot{\mathrm{x}}(s) = f(s,\mathrm{x}(s),\alpha(s)), s\in[t,T], \text{  and  } \mathrm{x}(t)=x,
    \label{eq:dynamics}
\end{align}
where $(t,x)$ are the initial time and state, $\mathrm{x}:[t,T]\times\R^n\times A\rightarrow\R^n$ is the state trajectory, $A(t)\subset \R^m$ is the control constraint, $\alpha\in\mathcal{A}(t)$ is the control signal, and we denote the set of measurable control signals 
\begin{align}
    \mathcal{A}(t) \coloneqq \{\alpha:[t,T]\rightarrow A ~|~ \|\alpha\|_{L^\infty(t,T)} <\infty\}.
    \label{eq:def_CtrlTraj}
\end{align}
We assume that $A$ is a compact subset in $\R^m$. 

With this dynamical system, we solve two optimal control problems.
\begin{problem} For given initial time and state $(t,x)$, solve 
\begin{align}
    \begin{split}
    & \vartheta_1 (t,x) \coloneqq \inf_{\alpha\in\mathcal{A}} \max_{\tau\in[t,T]} \int_t^\tau L(s,\mathrm{x}(s),\alpha(s))ds  \\
    &\quad\quad\quad\quad\quad\quad\quad\quad\quad\quad\quad\quad\quad\quad\quad  + g(\tau,\mathrm{x}(\tau)), 
    \label{eq:def_vartheta1}
    \end{split}
    \\
    & \text{subject to } \quad c(s,\mathrm{x}(s)) \leq 0, \quad s\in [t,T],
    \label{eq:def_vartheta_const1}
\end{align}
where $\mathrm{x}$ solves \eqref{eq:dynamics}.
\label{prob:1}
\end{problem}
\begin{problem} For given initial time and state $(t,x)$, solve 
\begin{align}
    \begin{split}
    & \vartheta_2 (t,x) \coloneqq \inf_{\alpha\in\mathcal{A}} \min_{\tau\in[t,T]} \int_t^\tau L(s,\mathrm{x}(s),\alpha(s))ds  \\
    &\quad\quad\quad\quad\quad\quad\quad\quad\quad\quad\quad\quad\quad\quad\quad  + g(\tau,\mathrm{x}(\tau)), 
    \label{eq:def_vartheta2}
    \end{split}
    \\
    & \text{subject to } \quad c(s,\mathrm{x}(s)) \leq 0, \quad s\in [t,\tau].
    \label{eq:def_vartheta_const2}
\end{align}
where $\mathrm{x}$ solves \eqref{eq:dynamics}.
\label{prob:2}
\end{problem}
For both problems, $L:[t,T]\times \R^n\times A \rightarrow \R$ is the stage cost, $g:\R\times\R^n\rightarrow\R$ is the terminal cost, $f:[t,T]\times\R^n\times A\rightarrow\R^n$ is the system dynamics, and $c:[t,T]\times\R^n\rightarrow\R$ is the state constraint.

In Problem \ref{prob:1}, we are concerned with finding a control signal that minimizes the maximum cost over the time horizon under the constraint. 
On the other hand, in Problem \ref{prob:2}, we determine a control signal that minimizes the minimum cost over the time horizon under the constraint.
The HJ equations for Problems \ref{prob:1} and \ref{prob:2} were presented in \cite{donggunTAC_part1} and \cite{donggun2020CDC}, respectively, and we review these in Section \ref{sec:HJeq}.

In this paper, we assume the following.
\begin{assumption}[Lipschitz continuity and compactness]
~
\begin{enumerate}
    \item The control set $A$ is compact and convex;
    
    \item $f:[0,T]\times\R^n\times A\rightarrow\R^n$  ($f=f(t,x,a)$) is Lipschitz continuous in $(t,x)$ for each $a\in A$;
    
    \item the stage cost $L:[0,T]\times\R^n\times \R^m\rightarrow\R$ ($L=L(t,x,a)$) is Lipschitz continuous in $(t,x)$ for each $a\in A$;
    
    \item for all $(t,x)\in[0,T]\times\R^n$, $\{f(t,x,a)~|~a\in A\}$ and $\{L(t,x,a)~|~a\in A\}$ are compact and convex;
    
    \item the terminal cost $g:[0,T]\times\R^n\rightarrow\R$ ($g=g(t,x)$) is Lipschitz continuous in $(t,x)$;
    
    \item the state constraint $c:[0,T]\times\R^n \rightarrow\R$ ($c=c(t,x)$) is Lipschitz continuous in $(t,x)$;
    
    \item the stage cost ($L$) and the terminal cost ($g$) are bounded below.
\end{enumerate}
\label{assum:BigAssum}
\end{assumption}
Assumption \ref{assum:BigAssum} guarantees the existence of the unique solution to the HJ equations, presented in Section \ref{sec:HJeq}.



\section{Hamilton-Jacobi Equations for Optimal Control Problems}
\label{sec:HJeq}

This section addresses the HJ equations for Problems \ref{prob:1} and \ref{prob:2} from the previous studies \cite{donggunTAC_part1,donggun2020CDC}.
For each Problem, two HJ equations are presented: one for time-varying stage cost, terminal cost, dynamics, and state constraint, the other for the time-invariant case: stage cost, terminal cost, dynamics, and state constraint are time-invariant.


\subsection{Hamilton-Jacobi equation for Problem \ref{prob:1}}
\label{sec:HJ_Prob1}

Instead of solving $\vartheta_1$ in \eqref{eq:def_vartheta1} subject to \eqref{eq:def_vartheta_const1}, we derive an HJ equation for a different value function $V_1:[0,T]\times\R^n\times\R$
\begin{align}
\begin{split}
    V_1&(t,x,z) \coloneqq    \inf_{\alpha\in\mathcal{A}(t)} \max\Big\{  \max_{s\in[t,T]}c(s,\mathrm{x}(s)),\\
    &\max_{\tau\in[t,T]} \int_t^\tau L(s,\mathrm{x}(s),\alpha(s))ds+g(\tau,\mathrm{x}(\tau))-z \Big\},
\end{split}
    \label{eq:ctrl_valueDef1}
\end{align}
where $\mathrm{x}$ solves \eqref{eq:dynamics} and $z$ is an auxiliary variable in $\R$.
The reason to solve \eqref{eq:ctrl_valueDef1} instead of \eqref{eq:def_vartheta1} subject to \eqref{eq:def_vartheta_const1} is that $\vartheta_1$ is not continuous in $(t,x)$-space since $\vartheta_1$ is infinity if there is no control signal to satisfy the state constraint.
Thus, viscosity theory cannot be applied. 
Note that the viscosity theory deals with a concept of weak solutions to nonlinear-first-order PDEs, and viscosity solutions are continuous and bounded.

On the other hand, $V_1$ is continuous in $(t,x,z)$-space, and $V_1$ encodes information of $\vartheta_1$ as stated in Lemma \ref{lemma:Equiv_twoCost1}.

\begin{lemma}[Equivalence of two value functions \cite{donggunTAC_part1}] Suppose Assumption \ref{assum:BigAssum} holds. For all $(t,x)\in[0,T]\times\R^n$, 
    \begin{align}
        \vartheta_1(t,x) = \min z \text{ subject to } V_1(t,x,z)\leq 0,
        \label{eq:lemma_equiv_twoCost1}
    \end{align}
    where $\vartheta_1$ and $V_1$ are defined in \eqref{eq:def_vartheta1} subject to \eqref{eq:def_vartheta_const1}, and in \eqref{eq:ctrl_valueDef1}, respectively.
    \label{lemma:Equiv_twoCost1}
\end{lemma}

Theorem \ref{thm:HJeq_prob1} presents a corresponding HJ equation for $V_1$.

\begin{theorem}[HJ equation for Problem \ref{prob:1} \cite{donggunTAC_part1}]
    Suppose Assumption \ref{assum:BigAssum} holds. 
    $V_1$ in \eqref{eq:ctrl_valueDef1} is the unique viscosity solution to the HJ equation:
    \begin{align}
    \begin{split}
        \max\big\{c(t,x)-&V_1(t,x,z), g(t,x)-z-V_1(t,x,z), \\ 
        &V_{1,t} - \bar{H} (t,x,z,D_x V_1,D_z V_1)\big\}=0
        \end{split}
        \label{eq:HJeq1}
    \end{align}
    in $(0,T)\times\R^n\times\R$, where $\bar{H}:[0,T]\times\R^n\times\R\times\R^n\times\R\rightarrow\R$
    \begin{align}
        \bar{H}(t,x,z,p,q) \coloneqq \max_{a\in A} - p\cdot f(t,x,a) + q L(t,x,a),
        \label{eq:def_Hambar}
    \end{align}
    and 
    \begin{align}
        V_1(T,x,z) = \max\{ c(T,x), g(T,x)-z \}
        \label{eq:V_terminalValue1}
    \end{align}
    on $\{t=T\}\times\R^n\times\R$. Denote $V_{1,t}=\frac{\partial V_1}{\partial t}$, $D_x V_1 = \frac{\partial V_1}{\partial x}$, and $D_z V_1 = \frac{\partial V_1}{\partial z}$.
    \label{thm:HJeq_prob1}
    
    For the time-invaraint case,
    $V_1$ is also the unique viscosity solution to the HJ equation:
    \begin{align}
    \begin{split}
        \max\big\{c(x)-&V_1 (t,x,z), \\
        &V_{1,t} - \bar{H}_1^\textnormal{TI} (x,z,D_x V_1,D_z V_1)\big\} =0
    \end{split}
    \label{eq:HJeq1_TI}
    \end{align}
    in $(0,T)\times\R^n\times\R$,
    where 
    \begin{align}
        \bar{H}^{\textnormal{TI}}_1(x,z,p,q) \coloneqq \min\big\{0,\bar{H} (x, z, p, q)\big\}
        \label{eq:Ham1_TI}
    \end{align}
    for $(x,z,p,q)\in\R^n\times\R\times\R^n\times\R$, 
    $\bar{H}$ is defined in \eqref{eq:def_Hambar} and the time dependency is omitted for the time-invariant case,
    and
    \begin{align}
        V_1 (T,x,z) = \max\{c(x),g(x)-z\}
    \end{align}
    on $\{t=T\}\times\R^n\times\R$.
\end{theorem}

This theorem states that, for the time-invariant case, the HJ equation \eqref{eq:HJeq1} is equivalent to \eqref{eq:HJeq1_TI}.


\subsection{Hamilton-Jacobi equation for Problem \ref{prob:2}}
\label{sec:HJ_Prob2}

The HJ analysis for Problem \ref{prob:2} is similar to that for Problem \ref{prob:1}.
With an auxiliary variable $z\in\R$, we define a new value function $V_2$ that combines the cost \eqref{eq:def_vartheta2} and the constraint \eqref{eq:def_vartheta_const2} of $\vartheta_2$: for $(t,x,z)\in[0,T]\times\R^n\times\R$,
\begin{align}
\begin{split}
    V_2(t,x,z) \coloneqq    \inf_{\alpha\in\mathcal{A}(t)} \min_{\tau\in[t,T]}\max\Big\{  \max_{s\in[t,\tau]}c(s,\mathrm{x}(s)),\\
    \int_t^\tau L(s,\mathrm{x}(s),\alpha(s))ds+g(\tau,\mathrm{x}(\tau))-z \Big\},
\end{split}
    \label{eq:ctrl_valueDef2}
\end{align}
where $\mathrm{x}$ solves \eqref{eq:dynamics}.
Note that $V_1$ in \eqref{eq:ctrl_valueDef1} and $V_2$ in \eqref{eq:ctrl_valueDef2} are different: $V_1$ contains the maximum operation over $\tau$, and $V_2$ has the minimum operation over $\tau$.

$\vartheta_2$ and $V_2$ satisfy the following lemma.
\begin{lemma}[Equivalence of two value functions \cite{donggun2020CDC}] Suppose Assumption \ref{assum:BigAssum} holds. For all $(t,x)\in[0,T]\times\R^n$, 
    \begin{align}
        \vartheta_2(t,x) = \min z \text{ subject to } V_2(t,x,z)\leq 0,
        \label{eq:lemma_equiv_twoCost2}
    \end{align}
    where $\vartheta_2$ and $V_2$ are defined in \eqref{eq:def_vartheta2} subject to \eqref{eq:def_vartheta_const2}, and in \eqref{eq:ctrl_valueDef2}.
    \label{lemma:Equiv_twoCost2}
\end{lemma}
Theorem \ref{thm:HJeq_prob2} presents an HJ equation for $V_2$.

\begin{theorem}[HJ equation for Problem \ref{prob:2} \cite{donggun2020CDC}]
    Suppose Assumption \ref{assum:BigAssum} holds. $V_2$ in \eqref{eq:ctrl_valueDef2} is the unique viscosity solution to the HJ equation:
    \begin{align}
        \max\Big\{c(t,x)- &V_2(t,x,z), \min\big\{g(t,x)-z-V_2(t,x,z), \notag\\
        &V_{2,t} - \bar{H}(t,x,z,D_x V_2,D_z V_2)\big\}\Big\}=0
    \label{eq:HJeq2}
    \end{align}
    in $(0,T)\times\R^n\times\R$, where 
    $\bar{H}$ is defined in \eqref{eq:def_Hambar},
    and 
    \begin{align}
        V_2(T,x,z) = \max\{ c(T,x), g(T,x)-z \}
        \label{eq:V_terminalValue2}
    \end{align}
    on $\{t=T\}\times\R^n\times\R$.
    
    For the time-invariant case,
    $V_2$ is also the unique viscosity solution to HJ equation:
    \begin{align}
    \begin{split}
        \max\big\{c(x)-&V_2 (t,x,z), \\
        &V_{2,t} - \bar{H}_2^\textnormal{TI} (x,z,D_x V_2,D_z V_2)\big\} =0
    \end{split}
    \label{eq:HJeq2_TI}
    \end{align}
    in $(0,T)\times\R^n\times\R$,
    where 
    \begin{align}
        \bar{H}_2^\textnormal{TI} (x,z,p,q) = \max\big\{0,\bar{H}(x,z,p,q)\big\}
        \label{eq:Ham2_TI}
    \end{align}
    for $(x,z,p,q)\in\R^n\times\R\times\R^n\times\R$,
    $\bar{H}$ is defined in \eqref{eq:def_Hambar},
    and
    \begin{align}
        V_2 (T,x,z) = \max\{c(x),g(x)-z\}
        \label{eq:V_TI_terminalValue2}
    \end{align}
    on $\{t=T\}\times\R^n\times\R$.
    \label{thm:HJeq_prob2}
\end{theorem}



The Hamiltonian $\bar{H}_2^\textnormal{TI}$ for the time-invariant Problem \ref{prob:2} in \eqref{eq:HJeq2_TI} is $\max\{0,\bar{H}\}$, which is convex in the costate space $(p,q)$ for each $(x,z)$.
On the other hand, the Hamiltonian $\bar{H}_1^\textnormal{TI}$ for the time-invariant Problem \ref{prob:1} in \eqref{eq:HJeq1_TI} is $\min\{0,\bar{H}\}$, which is not convex in the costate space $(p,q)$.

\subsection{Computational complexity of grid-based method}

For Problems \ref{prob:1} and \ref{prob:2}, the HJ equations in Theorems \ref{thm:HJeq_prob1} and \ref{thm:HJeq_prob2} can be numerically solved by grid-based methods, such as the level-set methods \cite{Osher02} and fast marching method \cite{sethian1996fast}.
These methods require spatial and temporal discretization, which leads to computational complexity exponential in the dimension of the state. 
Thus, it is intractable to utilize these grid-based methods for systems of state dimension beyond six.

\section{Lax formula for the state-constrained optimal control problem}
\label{sec:Laxformula}

In this section, we propose Lax formulae for Problems \ref{prob:1} and \ref{prob:2}.
In Section \ref{sec:convexity_valuefunction}, we first propose viscosity theory, which will be utilized to derive Lax formulae for these problems in the following subsections: Section \ref{sec:Laxformula_derivation1} is for Problem \ref{prob:1}; Section \ref{sec:Laxformula_derivation2} is for Problem \ref{prob:2}; and Section \ref{sec:Laxformula_derivation_TI_2} is for the time-invariant Problem \ref{prob:2}.
There is no Lax formula for the time-invariant Problem \ref{prob:1} since the corresponding Hamiltonian \eqref{eq:Ham1_TI} is non-convex in the costate.
This issue will be discussed more in Section \ref{sec:Laxformula_derivation_TI_2}.


\subsection{Sufficient condition under which two different PDEs have the same solution }
\label{sec:convexity_valuefunction}

This section proposes a general technique in viscosity theory to investigate the equivalence of two first-order PDEs.

Suppose $t\in[0,T]\subset\R$, $x\in\R^{n_x}$, $z\in\R^{n_z}$. 
Consider two first-order differential functions: for $i=1,2$, $F_i(t,x,z,u,r,p,q):[0,T]\times\R^{n_x}\times\R^{n_z}\times\R\times\R\times\R^{n_x}\times\R^{n_z}\rightarrow\R$, and suppose $X_i(t,x,z):[0,T]\times\R^{n_x}\times\R^{n_z}\rightarrow\R$ is the unique viscosity solution to 
\begin{align}
\begin{split}
    0 = F_i(t,x,z&,X_i(t,x,z), X_{i,t}(t,x,z),\\
    &D_x X_i(t,x,z),D_z X_i(t,x,z))
\end{split}
    \label{eq:FirstOrder_PDE}
\end{align}
in $(t,x,z)\in(0,T)\times\R^{n_x}\times\R^{n_z}$, and the terminal values for $X_1$ and $X_2$ are the same:
\begin{align}
    X_1(T,x,z)=X_2(T,x,z)\quad  \forall (x,z) \in \R^{n_x}\times\R^{n_z}.
    \label{eq:FirstOrder_TerminalValue}
\end{align}
We present conditions under which $X_1$ and $X_2$ are the same even though $F_1$ and $F_2$ are different.

Consider superdifferentials and subdifferentials of $X_i$ ($i=1,2$) with respect to $z$:
for each $(t,x,z)\in[0,T]\times\R^{n_x}\times\R^{n_z}$, $q\in D_z^+ X_i(t,x,z)$ ($i=1,2$) is a superdifferential with respect to $z$, if
\begin{align}
    D_z^+ X_i(t,x,z)\coloneqq\big\{q~|~ \limsup_{\bar{z}\rightarrow 0}  & \big(X_i(t,x,z+\bar{z}) - X_i(t,x,z) \notag\\
    &- q\cdot \bar{z} \big)/\lVert \bar{z} \rVert \leq 0 \},
    \label{eq:ineq_supergrad}
\end{align}
and $q\in D_z^- X_i(t,x,z)$ ($i=1,2$) is a subdifferential with respect to $z$, if 
\begin{align}
    D_z^- X_i(t,x,z)\coloneqq\big\{q~|~ \liminf_{\bar{z}\rightarrow 0}  & \big(X_i(t,x,z+\bar{z}) - X_i(t,x,z) \notag\\
    &- q\cdot \bar{z} \big)/\lVert \bar{z} \rVert \geq 0 \}.
    \label{eq:ineq_subgrad}
\end{align}



Theorem \ref{thm:ViscosityTheory} generalizes the condition for the equivalence of the two value functions ($X_1$ and $X_2$) even in the case where $F_1$ is different from $F_2$ for some $(t,x,z,u,r,p,q)$.
\begin{theorem}
    Suppose each of the two first-order differential equations in \eqref{eq:FirstOrder_PDE} with the terminal value \eqref{eq:FirstOrder_TerminalValue} for $i=1,2$ has the unique solution ($X_i$). 
    If, for all $(t,x,z,u,r,p)\in[0,T]\times\R^{n_x}\times\R^{n_z}\times\R\times\R\times\R^{n_x}$ 
    \begin{align}
    \begin{split}
        F_1 (t,x,z,u,r&,p,q) = F_2 (t,x,z,u,r,p,q)    \\
        &\forall q\in D_z^+ X_1(t,x,z) \cup D_z^- X_1 (t,x,z),
    \end{split}
        \label{eq:ViscosityThm_Condition}
    \end{align}
    then $X_1 \equiv X_2$.
    \label{thm:ViscosityTheory}
\end{theorem}
\textbf{Proof.} See Appendix \ref{appen:thm_ViscosityTheory}.

\subsection{Lax formula for Problem \ref{prob:1}}
\label{sec:Laxformula_derivation1}

In this subsection, we utilize the theory in Section \ref{sec:convexity_valuefunction} to derive and prove a Lax formula for Problem \ref{prob:1}, the solution to the HJ equation in Theorem \ref{thm:HJeq_prob1}.


We first investigate the superdifferentials and subdifferentials of $V_i$ ($i=1,2$) with respect to $z$.
\begin{lemma}[Convexity of the value function in $z$]
    For each $(t,x)\in [0,T] \times \R^n$, $V_1(t,x,\cdot)$ in \eqref{eq:ctrl_valueDef1} and $V_2(t,x,\cdot)$ in \eqref{eq:ctrl_valueDef2} are convex in $z\in\R$: for all $z_1,z_2\in\R$ and $\theta\in[0,1]$,
    \begin{align}
        V_i(t,x, \theta_1 z_1 + \theta_2 z_2) \leq \theta_1 V_i(t,x,z_1) + \theta_2 V_i(t,x,z_2)
        \label{eq:cvx_V_z}
    \end{align}
    for $i=1,2$.
    \label{lemma:cvx_V_z}
\end{lemma}
\textbf{Proof.} See Appendix \ref{appen:cvx_V_z}.
\begin{lemma}
    For all $(t,x,z)\in[0,T]\times\R^n\times\R$ and $i=1,2$, 
    \begin{align}
        D_z^- V_i(t,x,z) \subset [-1,0],
        \label{eq:lemma2_1}
    \end{align}
    and if $D_z^+ V_i(t,x,z)$ is not the empty set, there exists a unique $q\in\R$ such that
    \begin{align}
        D_z^+ V_i(t,x,z) = \{D_z V_i(t,x,z)\} \subset [-1,0].
        \label{eq:lemma2_2}
    \end{align}
    Note that $V_1$ and $V_2$ are defined in \eqref{eq:ctrl_valueDef1} and \eqref{eq:ctrl_valueDef2}, respectively.
    \label{lemma:bound_gradient}
\end{lemma}
\textbf{Proof.} See Appendix \ref{appen:lemma_bound_gradient}.



Define, 
for initial time and state $(t,x)\in[0,T]\times\R^n$,
\begin{align}
    \begin{split}
    & \varphi_1 (t,x) \coloneqq \inf_{\beta} \max_{\tau\in[t,T]} \int_t^\tau H^*(s,\mathrm{x}(s),\beta(s))ds  \\
    &\quad\quad\quad\quad\quad\quad\quad\quad\quad\quad\quad\quad\quad\quad\quad  + g(\tau,\mathrm{x}(\tau)), \label{eq:def_varphi1}
    \end{split}
    \\
    & \text{subject to }
    \begin{cases}
        \dot{\mathrm{x}}(s) = -\beta(s), & s\in [t,T], \\
        \mathrm{x}(t) = x,\\
        c(s,\mathrm{x}(s)) \leq 0, & s\in [t,T],
    \end{cases}
    \label{eq:def_varphi_const1}
\end{align}
where $H:[0,T]\times\R^n\times\R^n\rightarrow\R$
\begin{align}
    H(s,x,p) \coloneqq  \bar{H}(s,x,z,p,-1),
    \label{eq:def_Ham}
\end{align}
and $H^*:[0,T]\times\R^n\times\R^n\rightarrow\R$
\begin{align}
    H^*(s,x,b) \coloneqq \max_{p\in\R^n} p\cdot b - H(s,x,p).
    \label{eq:def_ConvConj_Ham}
\end{align}
$H^*$ is the Legendre-Fenchel transformation (the convex conjugate) of $H$ with respect to $p\in\R^n$, and the state constraint, $c(s,\mathrm{x}(s))\leq 0$, is satisfied in $[t,T]$ but not $[t,\tau]$.
Combining \eqref{eq:def_varphi1} and \eqref{eq:def_varphi_const1}, define a value function, $W_1:[0,T]\times\R^n\times\R\rightarrow\R$:
\begin{align}
\begin{split}
    W_1&(t,x,z) = \inf_{\beta} \max\Big\{\max_{s\in[t,T] }c(s,\mathrm{x}(s)) , \\
    &\max_{\tau\in[t,T]} \int_t^\tau H^*(s,\mathrm{x}(s),\beta(s))ds + g(\tau,\mathrm{x}(\tau))-z \Big\}.
    \label{eq:ctrl_W1}
\end{split}
\end{align}
By Theorem \ref{thm:HJeq_prob1}, $W_1$ is the unique viscosity solution to 
\begin{align}
\begin{split}
    \max\big\{c(t,x)&-W_1(t,x,z),g(t,x)-z-W_1(t,x,z), \\ 
    &W_{1,t} - \bar{H}_W(t,x,z,D_x W_1,D_z W_1)\big\}=0
    \label{eq:HJeqW1}
\end{split}
\end{align}
in $(0,T)\times\R^n\times\R$, where $\bar{H}_W:[0,T]\times\R^n\times\R\times\R^n\times\R\rightarrow\R$
\begin{align}
    \bar{H}_W(t,x,z,p,q) \coloneqq \max_{b} p\cdot b +q H^*(t,x,b),
    \label{eq:def_Hambar1}
\end{align}
$H^*$ is defined in \eqref{eq:def_ConvConj_Ham},
and 
\begin{align}
    W_1(T,x,z) = \max\{ c(T,x), g(T,x)-z \}
    \label{eq:W_terminalValue1}
\end{align}
on $\{t=T\}\times\R^n\times\R$.

In the rest of this subsection, we build a mathematical background to utilize Theorem \ref{thm:ViscosityTheory} in order to prove that $\vartheta_1$ (\eqref{eq:def_vartheta1} subject to \eqref{eq:def_vartheta_const1}) and $\varphi_1$ (\eqref{eq:def_varphi1} subject to \eqref{eq:def_varphi_const1}) are the same.


For $s\in[t,T]$, the control $\beta(s)$ is constrained in the domain of $H^*(s,\mathrm{x}(s),\cdot)$,
\begin{align}
    \text{Dom}(H^*(s,\mathrm{x}(s),\cdot)) = \{ b~|~ H^*(s,\mathrm{x}(s),b)<\infty \}.
    \label{eq:def_Dom}
\end{align}
which is specified in Lemma \ref{lemma:beta_const}.
\begin{lemma}[Domain of $H^*$ \cite{donggunTAC_submitted}]
    For all $(t,x)\in[0,T]\times\R^n$,
    \begin{align}
        \text{Dom}(H^*(t,x,\cdot)) = \text{co}(B(t,x)),
    \end{align}
    where $H^*$ is defined in \eqref{eq:def_ConvConj_Ham}, for $(t,x)\in[0,T]\times\R^n$,
    \begin{align}
        B(t, x) \coloneqq \{-f(t,x,a)~|~a\in A\},
        \label{eq:ctrlbound_B}
    \end{align}
    and $\text{co}(B(t,x))$ is the convex hull of the set $B(t,x)$.
    \label{lemma:beta_const}
\end{lemma}
Although the control constraint for $\beta(s)$ is omitted in \eqref{eq:def_varphi_const1}, $\beta(s)$ is constrained in $\text{co}(B(s,\mathrm{x}(s)))$ (the domain of $H^*(t,x,\cdot)$).
By combining Lemma \ref{lemma:beta_const} and the definitions of Hamiltonians ($\bar{H}$ in \eqref{eq:def_Hambar} and $\bar{H}_W$ in \eqref{eq:def_Hambar1}), we derive the following lemma.
\begin{lemma}
    For $(t,x,z,p,q)\in[0,T]\times\R^n\times\R\times\R^n\times\R$,
    \begin{align}
        \bar{H}(t,x,z,p,q) = \bar{H}_W(t,x,z,p,q)\quad \text{ if } q\leq0,
        \label{eq:Hamiltonian_comp}
    \end{align}
    where $\bar{H}$ and $\bar{H}_W$ are defined in \eqref{eq:def_Hambar} and \eqref{eq:def_Hambar1}, respectively.
    \label{lemma:Hamiltonian_comp}
\end{lemma}
\textbf{Proof.} See Appendix \ref{appen:lemma_Hamiltonian_comp}.

For $q>0$, \eqref{eq:Hamiltonian_comp} does not generally hold.
Therefore, the two HJ equations in \eqref{eq:HJeq1} and \eqref{eq:HJeqW1} are different in general. 
However, as proved in Lemma \ref{lemma:bound_gradient}, the subdifferentials and superdifferentials of $V_1(t,x,z)$ with respect to $z$ are less than or equal to 0 for all $(t,x,z)$.
Thus, by combining Theorem \ref{thm:ViscosityTheory}, Lemma \ref{lemma:bound_gradient} for $i=1$, and Lemma \ref{lemma:Hamiltonian_comp}, we prove $V_1\equiv W_1$ as in Theorem \ref{thm:Laxformula1}.
This implies that $\vartheta_1$ (\eqref{eq:def_vartheta1} subject to \eqref{eq:def_vartheta_const1}) and $\varphi_1$ (\eqref{eq:def_varphi1} subject to \eqref{eq:def_varphi_const1}) are the same by Lemma \ref{lemma:Equiv_twoCost1}.
This paper calls $\varphi_1$ the Lax formula for Problem \ref{prob:1}. 

\begin{theorem}\textnormal{\textbf{(Lax formula for Problem \ref{prob:1})}}
    For all $(t,x)\in[0,T]\times\R^n$, $W_1$ \eqref{eq:ctrl_W1} is the unique viscosity solution to the HJ equation  \eqref{eq:HJeq1} and \eqref{eq:V_terminalValue1} in Theorem \ref{thm:HJeq_prob1}, i.e.
    \begin{align}
        V_1 \equiv W_1.
    \end{align}
    Also, 
    \begin{align}
        \vartheta_1(t,x) = \varphi_1(t,x),
    \end{align}
    where $\vartheta_1$ is \eqref{eq:def_vartheta1} subject to \eqref{eq:def_vartheta_const1}, and $\varphi_1$ is \eqref{eq:def_varphi1} subject to \eqref{eq:def_varphi_const1}.
    \label{thm:Laxformula1}
\end{theorem}

\begin{remark}
    This paper calls $\varphi_1$ in \eqref{eq:def_varphi1} subject to \eqref{eq:def_varphi_const1} the Lax formula for Problem \ref{prob:1}, which provides the optimal value $\vartheta_1$ of Problem \ref{prob:1}.
\end{remark}

\subsection{Lax formula for Problem \ref{prob:2}}
\label{sec:Laxformula_derivation2}

This subsection presents a Lax formula for Problem \ref{prob:2}.
We use an analogous derivation presented in Section \ref{sec:Laxformula_derivation1}. 

Define, for initial time and state $(t,x)\in[0,T]\times\R^n$, 
\begin{align}
    \begin{split}
    & \varphi_2 (t,x) \coloneqq \inf_{\beta} \min_{\tau\in[t,T]} \int_t^\tau H^*(s,\mathrm{x}(s),\beta(s))ds  \\
    &\quad\quad\quad\quad\quad\quad\quad\quad\quad\quad\quad\quad\quad\quad\quad  + g(\tau,\mathrm{x}(\tau)), \label{eq:def_varphi2}
    \end{split}
    \\
    & \text{subject to }
    \begin{cases}
        \dot{\mathrm{x}}(s) = -\beta(s), & s\in [t,T], \\
        \mathrm{x}(t) = x,\\
        c(s,\mathrm{x}(s)) \leq 0, & s\in [t,\tau],
    \end{cases}
    \label{eq:def_varphi_const2}
\end{align}
where $H^*$ is defined in \eqref{eq:def_ConvConj_Ham}, and $\beta(s)$ is constrained in the domain of $H^*(s,\mathrm{x}(s),\cdot)$. 
We call $\varphi_2$ the Lax formula for Problem \ref{prob:2}.


Following the same argument for the proof of Theorem \ref{thm:Laxformula1}, Theorem \ref{thm:Laxformula2} is derived.


\begin{theorem}\textnormal{\textbf{(Lax formula for Problem \ref{prob:2})}}
    For all $(t,x)\in[0,T]\times\R^n$,
    \begin{align}
        \vartheta_2(t,x) = \varphi_2(t,x),
    \end{align}
    where $\vartheta_2$ is \eqref{eq:def_vartheta2} subject to \eqref{eq:def_vartheta_const2}, and $\varphi_2$ is \eqref{eq:def_varphi2} subject to \eqref{eq:def_varphi_const2}.
    \label{thm:Laxformula2}
\end{theorem}

\begin{remark}
    This paper calls $\varphi_2$ in \eqref{eq:def_varphi2} subject to \eqref{eq:def_varphi_const2} the Lax formula for Problem \ref{prob:2}, which provides the optimal value $\vartheta_2$ of Problem \ref{prob:2}.
\end{remark}

\subsection{Lax formula for the time-invariant Problem \ref{prob:2}}
\label{sec:Laxformula_derivation_TI_2}

For the time-invariant case, there is no time dependency on the stage and terminal costs, dynamics, and state constraint.
Theorems \ref{thm:HJeq_prob1} and \ref{thm:HJeq_prob2} present HJ equations for Problem \ref{prob:1} and \ref{prob:2} for the time-invariant case, respectively.
On the derivation of Lax formulae, the Hamiltonian has to be convex in the costate, unless Lemma \ref{lemma:Hamiltonian_comp} does not hold.
The Hamiltonian for the time-invariant Problem \ref{prob:1}, $\bar{H}_1^\textnormal{TI}$ in \eqref{eq:Ham1_TI}, is generally non-convex, thus there is no corresponding Lax formula.
On the other hand, the Hamiltonian for the time-invariant Problem \ref{prob:2}, $\bar{H}_2^\textnormal{TI}$ in \eqref{eq:Ham2_TI}, is convex, thus the corresponding Lax formula can be derived.

Define, for the initial time and state $(t,x)\in[0,T]\times\R^n$,
\begin{align}
    \varphi_2^\textnormal{TI}(t,x)&\coloneqq \inf_\beta \int_t^T H_2^\textnormal{TI*}(\mathrm{x}(s),\beta(s))ds+g(\mathrm{x}(T)),\label{eq:def_varphi2_TI}\\
    &\text{subject to }\begin{cases}\dot{\mathrm{x}}(s)=-\beta(s),&s\in[t,T],\\\mathrm{x}(t)=x,\\c(\mathrm{x}(s))\leq 0, & s\in[t,T], \end{cases}
    \label{eq:def_varphi2_TI_const}
\end{align}
where $H_2^\textnormal{TI}:\R^n\times\R^n\rightarrow\R$
\begin{align}
    H_2^\textnormal{TI}(x,p)\coloneqq \bar{H}_2^\textnormal{TI}(x,z,p,-1),
    \label{eq:Ham2_TI_x}
\end{align}
where $\bar{H}_2^\textnormal{TI}$ is defined in \eqref{eq:Ham2_TI}, and $H_2^{\textnormal{TI}*}:\R^n\times\R^n\rightarrow\R$
\begin{align}
    H_2^\textnormal{TI*}(x,b)\coloneqq \max_p p\cdot b - H_2^\textnormal{TI} (x,p).
    \label{eq:Ham_star2_TI_x}
\end{align}
Note that 
\begin{align}
    H_2^\textnormal{TI}(x,p) = \max\{0,H(x,p)\},
    \label{eq:Ham_Ham2TI}
\end{align}
where $H$ is defined in \eqref{eq:def_Ham},
and $\bar{H}_2^\textnormal{TI}$ in \eqref{eq:Ham2_TI_x} has no $z$-dependency, thus $H_2^\textnormal{TI}$ does not depend on $z$.

In comparison to the Lax formula for Problem \ref{prob:2}, the Lax formula for the time-invariant Problem \ref{prob:2} only contains the minimum operation over the control signal $\beta$, on the other hand, the Lax formula for Problem \ref{prob:2} contains the two minimum operations over the control signal $\beta$ and the terminal time $\tau$.

The new control signal $\beta$ is constrained to the domain of the stage cost $H_2^{\textnormal{TI}*}$, presented in Lemma \ref{lemma:domain_ctrl_TI_Prob2}.

\begin{lemma}[Domain of $H_2^{\textnormal{TI}*}$] 
    For all $x\in \R^n$,
    \begin{align}
        \text{Dom}(H_2^{\textnormal{TI}*}(x,\cdot)) = \text{co}(\{0\}\cup B(x)),
    \end{align}
    where $H_2^{\textnormal{TI}*}$ is defined in \eqref{eq:Ham_star2_TI_x}, 
    \begin{align}
        B(x) \coloneqq \{-f(x,a)~|~a\in A\},
        \label{eq:ctrlBound_B_TI}
    \end{align}
    which does not have the time dependency in comparison to $B(t,x)$ in \eqref{eq:ctrlbound_B}, and $\text{co}$ refers to the convex hull operation.
    \label{lemma:domain_ctrl_TI_Prob2}
\end{lemma}
\textbf{Proof.} See Appendix \ref{appen:lemma_domain_ctrl_TI_Prob2}.

We will prove that $\varphi_2^\textnormal{TI}$ are the same as $\vartheta_2$ for the time-invariant case by the viscosity theory in Section \ref{sec:convexity_valuefunction}.
Corresponding to $\varphi_2^\textnormal{TI}$, define a value function $W_2^\textnormal{TI}:[0,T]\times\R^n\times\R\rightarrow\R$:
\begin{align}
\begin{split}
    W_2^\textnormal{TI}(t&,x,z) = \inf_\beta\Big\{ \max_{s\in[t,T]}c(\mathrm{x}(s)),\\
    &\int_t^T H_2^{\textnormal{TI}*}(s,\mathrm{x}(s),\beta(s))ds+g(\mathrm{x}(T))-z\Big\}.
\end{split}
\label{eq:W_2^TI}
\end{align}
Theorem 3.1 in \cite{altarovici2013general} states the relationship between $W_2^\textnormal{TI}$ and $\varphi_2^\textnormal{TI}$ as below.
\begin{align}
    \varphi_2^\textnormal{TI} =\min z \text{  subject to  }W_2^\textnormal{TI} \leq 0
    \label{eq:equiv_phi2_TI_W2_TI}
\end{align}
Also, Proposition 3.4 in \cite{altarovici2013general} shows that $W_2^\textnormal{TI}$ is the unique viscosity solution to 
\begin{align}
\begin{split}
    \max\Big\{ &c(x)- W_{2}^\textnormal{TI}, \\
    &W_{2,t}^\textnormal{TI} - \bar{H}_W^\textnormal{TI}(x,z,D_x W_{2}^\textnormal{TI}, D_z W_{2}^\textnormal{TI})\Big\}=0
\end{split}
\label{eq:HJeq_W2_TI}
\end{align}
in $(0,T)\times\R^n\times\R$, where $\bar{H}_W^\textnormal{TI}:\R^n\times\R\times\R^n\times\R\rightarrow\R$
\begin{align}
    \bar{H}_W^\textnormal{TI}(x,z,p,q)\coloneqq \max_b p\cdot b +q H_2^{\textnormal{TI}*}(x,b),
    \label{eq:Ham_W_prob2_TI}
\end{align}
where $H_2^{\textnormal{TI}*}$ is defined in \eqref{eq:Ham_star2_TI_x}.

Lemma \ref{lemma:Hamiltonian_comp_TI} states that the Hamiltonian for $W_2^\textnormal{TI}$ ($\bar{H}_W^\textnormal{TI}$) in \eqref{eq:Ham_W_prob2_TI} is equivalent to the Hamiltonian for the time-invariant $V_2$ ($\bar{H}_2^\textnormal{TI}$) in \eqref{eq:Ham2_TI}. 

\begin{lemma}
    For $(x,z,p,q)\in\R^n\times\R\times\R^n\times\R$,
    \begin{align}
        \bar{H}_2^\textnormal{TI}(x,z,p,q) = \bar{H}_W^\textnormal{TI}(x,z,p,q)\quad \text{if }q\leq 0,
    \end{align}
    where $\bar{H}_2^\textnormal{TI}$ and $\bar{H}_W^\textnormal{TI}$ are defined in \eqref{eq:Ham2_TI} and \eqref{eq:Ham_W_prob2_TI}, respectively.
    \label{lemma:Hamiltonian_comp_TI}
\end{lemma}
\textbf{Proof.} See Appendix \ref{appen:lemma_Hamiltonian_comp_TI}.

Consider the two HJ equations \eqref{eq:HJeq2_TI} and \eqref{eq:HJeq_W2_TI}. By combining Lemma \ref{lemma:bound_gradient} for $i=2$, Lemma \ref{lemma:Hamiltonian_comp_TI}, and Theorem \ref{thm:ViscosityTheory}, we prove that $V_2$ in \eqref{eq:ctrl_valueDef2} and $W_2^\textnormal{TI}$ in \eqref{eq:W_2^TI} are the same viscosity solution to the HJ equation \eqref{eq:HJeq_W2_TI}, as presented in Theorem \ref{thm:Laxformula_Prob2_TI}.
Also, three value functions  ($\varphi_2^\textnormal{TI}$ in \eqref{eq:def_varphi2_TI} subject to \eqref{eq:def_varphi2_TI_const}, $\vartheta_2$ in \eqref{eq:def_vartheta2} subject to \eqref{eq:def_vartheta_const2}, and $\varphi_2$ in \eqref{eq:def_varphi2} subject to \eqref{eq:def_varphi_const2}) are the same by \eqref{eq:equiv_phi2_TI_W2_TI}, Lemma \ref{lemma:Equiv_twoCost2} for $i=2$, and Theorem \ref{thm:Laxformula2}.

\begin{theorem}\textnormal{\textbf{(Lax formula for Problem \ref{prob:2} (time-invariant version)) }}
    Consider Problem \ref{prob:2} for the time-invariant case: 
    $\vartheta_2$ in \eqref{eq:def_vartheta2} subject to \eqref{eq:def_vartheta_const2}.
For all $(t,x)\in[0,T]\times\R^n$, $W_2^\textnormal{TI}$ in \eqref{eq:W_2^TI} is the unique viscosity solution to the HJ equation in \eqref{eq:HJeq2_TI} and \eqref{eq:V_TI_terminalValue2}, i.e.,
\begin{align}
    V_2(t,x,z) = W_2(t,x,z) = W_2^\textnormal{TI}(t,x,z).
    \label{eq:valueSame_Prob2_TI}
\end{align}
In addition,
\begin{align}
    \vartheta_2(t,x) =\varphi_2(t,x)= \varphi_2^\textnormal{TI}(t,x),
\end{align}
where $\vartheta_2$ is \eqref{eq:def_vartheta2} subject to \eqref{eq:def_vartheta_const2}, $\varphi_2$ is \eqref{eq:def_varphi2} subject to \eqref{eq:def_varphi_const2}, and $\varphi_2^\textnormal{TI}$ is \eqref{eq:def_varphi2_TI} subject to \eqref{eq:def_varphi2_TI_const}.
\label{thm:Laxformula_Prob2_TI}
\end{theorem}
This paper calls $\varphi_2^\textnormal{TI}$ in \eqref{eq:def_varphi2_TI} subject to \eqref{eq:def_varphi2_TI_const} the Lax formula for the time-invariant Problem \ref{prob:2}. 


\section{Convexity analysis}
\label{sec:cvxAnalysis}

From this section, we set the initial time $t$ to 0 in Problems \ref{prob:1} and \ref{prob:2}.
Without spatial discretization, we temporally discretize problems and utilize gradient-based methods, which provide the optimal cost and a optimal control signal if the discretized problem is convex.
Thus, in this section, we present convexity analysis for the Problems \ref{prob:1} and \ref{prob:2} ($\vartheta_1$, $\vartheta_2$) and the corresponding Lax formulae ($\varphi_1$, $\varphi_2$, $\varphi_2^\textnormal{TI}$).
For the temporal discretization, the first-order forward Euler method is chosen.
Note that the convexity analysis in this section is invariant even if other temporal discretization methods are utilized, such as backward Euler method, Crank Nicolson method, and a family of Runge-Kutta methods.

\subsection{Convexity analysis for Problem \ref{prob:1} and the corresponding Lax formula}
\label{sec:Cvx_Analysis_Prob1}

This subsection presents a convexity analysis for the temporally discretized Problem \ref{prob:1}$ (\vartheta_1$ in \eqref{eq:def_vartheta1} subject to \eqref{eq:def_vartheta_const1}) and corresponding Lax formula ($\varphi_1$ in \eqref{eq:def_varphi1} subject to \eqref{eq:def_varphi_const1}).

The temporally discretized problem can be found using a temporal discretization $\{t_0=0,...,t_K=T\}$.
\begin{align}
    \begin{split}
    & \quad \vartheta_1 (0,x) \simeq \min_{\mathrm{x}[\cdot],\alpha[\cdot]} \max_{k'\in\{0,...,K\}} \sum_{k=0}^{k'} L(t_k,\mathrm{x}[k],\alpha[k])\Delta_k  \\
    &\quad\quad\quad\quad\quad\quad\quad\quad\quad\quad\quad\quad + g(t_{k'},\mathrm{x}[k']), 
    \label{eq:def_vartheta1_disc}
    \end{split}
    \\
    & \text{subject to } \begin{cases} \mathrm{x}[k+1]-\mathrm{x}[k] = \Delta_k f(t_k,\mathrm{x}[k],\alpha[k]),  \\ 
    \quad\quad\quad\quad\quad\quad\quad\quad\quad\quad\phantom{''} k\in\{0,...,K-1\},\\
    \alpha[k]\in A, \quad\quad\quad\quad\quad\quad\phantom{'} k\in\{0,...,K-1\},\\ 
    \mathrm{x}[0]=x,\\
    c(t_k,\mathrm{x}[k]) \leq 0, \quad\quad\quad\quad k\in \{0,...,K\}, 
    \end{cases}
    \label{eq:def_vartheta_const1_disc}
\end{align}
where $\Delta_k\coloneqq t_{k+1}-t_k$. In this paper, we use $\mathrm{x}$ to denote both a state trajectory and sequence. Note that $\mathrm{x}(s)$ is a state at time $s$, and $\mathrm{x}[k]$ is a $k$-th state.
This notation rule is also applied for $\alpha$.

The following lemma provides convexity conditions for the temporal discretized Problem \ref{prob:1} (\eqref{eq:def_vartheta1_disc}) subject to \eqref{eq:def_vartheta_const1_disc}) in $(\mathrm{x}[\cdot],\alpha[\cdot])$-space.

\begin{cvxCondition}\textnormal{\textbf{(Convexity conditions for the temporally discretized Problem \ref{prob:1} ($\vartheta_1$))}}
    \begin{enumerate}
        \item $L(t,x,a)$ is convex in $(x,a)$ for all $t\in[0,T]$,
        \item $g(t,x)$ is convex in $x$ for all $t\in[0,T]$,
        \item $c(t,x)$ is convex in $x$ for all $t\in[0,T]$,
        \item $f(t,x,a)$ is affine in $(x,a)$ for all $t\in[0,T]$,
        \item $A$ is convex.
    \end{enumerate}
\label{cond:cond1_vartheta1}
\end{cvxCondition}

\begin{lemma}
    Suppose Condition \ref{cond:cond1_vartheta1} holds. 
    For a fixed initial state $x\in\R^n$, 
    the temporally discretized Problem \ref{prob:1} ($\vartheta_1$ in \eqref{eq:def_vartheta1_disc} subject to \eqref{eq:def_vartheta_const1_disc}) is convex in $(\mathrm{x}[\cdot],\alpha[\cdot])$. 
    \label{lemma:cvx_vartheta1}
\end{lemma}
\textbf{Proof.} See Appendix \ref{appen:lemma_conv_vartheta1_cont}.

In this convexity analysis, we consider only $(\mathrm{x},\alpha)$ (or $(\mathrm{x}[\cdot],\alpha[\cdot])$) as variables to be determined, but not $\tau$ (or $k'$). 
The reason of choosing these decision variables is to avoid solving a minimax problem and to utilize the fact that the pointwise maximum of a family of convex functions is convex.

The temporally discretized Lax formula for Problem \ref{prob:1} is
\begin{align}
    \begin{split}
    & \varphi_1 (0,x) \simeq \min_{\mathrm{x}[\cdot],\beta[\cdot]} \max_{k'\in\{0,...,K\}} \sum_{k=0}^{k'} H^*(t_k,\mathrm{x}[k],\beta[k])\Delta_k  \\
    &\quad\quad\quad\quad\quad\quad\quad\quad\quad\quad\quad\quad + g(t_{k'},\mathrm{x}[k']), 
    \label{eq:def_varphi1_disc1}
    \end{split}
    \\
    & \text{subject to } \begin{cases} \mathrm{x}[k+1]-\mathrm{x}[k] = -\Delta_k \beta[k],  \\ 
    \quad\quad\quad\quad\quad\quad\quad\quad\quad\quad\phantom{''} k\in\{0,...,K-1\},\\
    \beta[k]\in \text{co}(B(t_k,\mathrm{x}[k])), \\
    \quad\quad\quad\quad\quad\quad\quad\quad\quad\quad\phantom{''} k\in\{0,...,K-1\},\\
    \mathrm{x}[0]=x,\\
    c(t_k,\mathrm{x}[k]) \leq 0, \quad\quad\quad\quad k\in \{0,...,K\}. \\  \end{cases}
    \label{eq:def_varphi_const1_disc1}
\end{align}
Lemma \ref{lemma:cvx_vartheta1} implies that if $H^*(s,x,b)$ and $\{(x,b)~|~b\in\text{co}(B(s,x))\}$ are convex in $(x,b)$ for each $s$, the temporally discretized Lax formula (\eqref{eq:def_varphi1_disc1} subject to \eqref{eq:def_varphi_const1_disc1}) is convex.
In one of the authors' paper \cite{donggunTAC_submitted}, sufficient conditions for convex $H^*(s,x,b)$ and $\{(x,b)~|~b\in\text{co}(B(s,x))\}$ are presented: 1) $L(s,x,a)=L^x(s,x)+L^a(s,a)$, and $L^x$ is convex in $x$ for all $s\in[0,T]$, 2) $f(s,x,a)=M(s)x+L^a(s,a)$ for some matrix $M(s)\in\R^n\times\R^n$.
Corollary \ref{corollary:cvx_varphi1} summarizes the convexity conditions for the temporally discretized Lax formula for Problem \ref{prob:1}.

\begin{cvxCondition} \textnormal{\textbf{(Convexity condition for the temporally discretized $\varphi_1$)}}
    \begin{enumerate}
        \item $L(s,x,a)=L^x(s,x)+L^a(s,a)$ for some $L^x$ and $L^a$, and $L^x$ is convex in $x$ for all $s\in[0,T]$,
        \item $g(s,x)$ is convex in $x$ for all $s\in[0,T]$,
        \item $c(s,x)$ is convex in $x$ for all $s\in[0,T]$,
        \item $f(s,x,a)=M(s)x+f^a(s,a)$ for some $M$ and $f^a$ for all $s\in[0,T]$.
    \end{enumerate}
\label{cond:cond1_varphi1}
\end{cvxCondition}
\begin{corollary}
    Suppose Condition \ref{cond:cond1_varphi1} holds. the Lax formula for Problem \ref{prob:1} ($\varphi_1$) in \eqref{eq:def_varphi1} subject to \eqref{eq:def_varphi_const1} is convex in $(\mathrm{x},\beta)$.
    Also, the temporally discretized Lax formula for Problem \ref{prob:1} (\eqref{eq:def_varphi1_disc1} subject to \eqref{eq:def_varphi_const1_disc1}) is convex in $(\mathrm{x}[\cdot],\beta[\cdot])$.
\label{corollary:cvx_varphi1}
\end{corollary}


\subsection{Convexity analysis for Problem \ref{prob:2} and the corresponding Lax formula}
\label{sec:Cvx_Analysis_Prob2}

This section presents convexity analysis for the temporal discretization of Problem \ref{prob:2}$ (\vartheta_2$ in  \eqref{eq:def_vartheta2} subject to \eqref{eq:def_vartheta_const2}), Lax formula for Problem \ref{prob:2} ($\varphi_2$ in \eqref{eq:def_varphi2} subject to \eqref{eq:def_varphi_const2}), and Lax formula for the time-invariant Problem \ref{prob:2} ($\varphi_2^\textnormal{TI}$ in  \eqref{eq:def_varphi2_TI} subject to \eqref{eq:def_varphi2_TI_const}).

The temporally discretized Problem \ref{prob:2} on $\{t_0=0,...,t_K=T\}$ is 
\begin{align}
    & \quad \vartheta_2 (0,x) \simeq \min_{\mathrm{x}[\cdot],\alpha[\cdot],k'\in\{0,...,K\}} \sum_{k=0}^{ k'} L(t_k,\mathrm{x}[k],\alpha[k])\Delta_k  \notag\\
    &\quad\quad\quad\quad\quad\quad\quad\quad\quad\quad\quad\quad + g(t_{k'},\mathrm{x}[k']), 
    \label{eq:def_vartheta2_disc2}
    \\
    & \text{subject to } \begin{cases} \mathrm{x}[k+1]-\mathrm{x}[k] = \Delta_k f(t_k,\mathrm{x}[k],\alpha[k]),  \\ 
    \quad\quad\quad\quad\quad\quad\quad\quad\quad\quad\phantom{''} k\in\{0,...,K-1\},\\
    \alpha[k]\in A, \quad\quad\quad\quad\quad\quad\phantom{'} k\in\{0,...,K-1\},\\
    \mathrm{x}[0]=x,\\
    c(t_k,\mathrm{x}[k]) \leq 0, \quad\quad\quad\quad k\in \{0,...,k' \}. \\  \end{cases}
    \label{eq:def_vartheta_const2_disc2}
\end{align}
In Section \ref{sec:Cvx_Analysis_Prob1}, we analyze the convexity conditions for the temporally discretized Problem \ref{prob:1} in $(\mathrm{x}[\cdot],\alpha[\cdot])$-space, on the other hand, in this section, the convexity conditions for Problem \ref{prob:2} is analyzed in $(\mathrm{x}[\cdot],\alpha[\cdot],k')$-space.
This is because Problem \ref{prob:2} is minimization problem over $(\mathrm{x}[\cdot],\alpha[\cdot],k')$-space.
The result of the convexity analysis is summarized in the fourth column of Table \ref{tab:comparisonCvxCond}.

Using the same temporal discretization $\{t_0=0,...,t_K=T\}$, the Lax formula for Problem \ref{prob:2} (\eqref{eq:def_varphi2} subject to \eqref{eq:def_varphi_const2}) is discretized to
\begin{align}
    \begin{split}
    & \varphi_2 (0,x) \simeq \min_{\mathrm{x}[\cdot],\beta[\cdot],k'\in\{0,...,K\}} \sum_{k=0}^{k'} H^*(s,\mathrm{x}[k],\beta[k])\Delta_k  \\
    &\quad\quad\quad\quad\quad\quad\quad\quad\quad\quad\quad\quad + g(t_{k'},\mathrm{x}[k']), 
    \label{eq:def_varphi2_disc}
    \end{split}
    \\
    & \text{subject to } \begin{cases} \mathrm{x}[k+1]-\mathrm{x}[k] = -\Delta_k \beta[k],  \\ 
    \quad\quad\quad\quad\quad\quad\quad\quad\quad\quad\quad\phantom{'''} k\in\{0,...,K-1\},\\
    \beta[k]\in \text{co}(B(t_k,\mathrm{x}[k])), \quad  k\in\{0,...,K-1\},\\
    \mathrm{x}[0]=x,\\
    c(t_k,\mathrm{x}[k]) \leq 0, \quad\quad\quad\quad\quad\phantom{''} k\in \{0,...,k'\}. \\  \end{cases}
    \label{eq:def_varphi_const2_disc}
\end{align}
Note that, \cite{donggunTAC_submitted} presents that $L=0$ implies $H^*=0$, and the control constraint ($\text{co}(B(t_k,\mathrm{x}[k]))$) is convex if $f(s,x,a)=M(s)x + f^a (s,a)$ for some $M$ and $f^a$ for all $s\in[0,T]$.
The result of the convexity analysis is summarized at the fifth column of Table \ref{tab:comparisonCvxCond}.

\subsection{Convexity analysis for the time-invariant Problem \ref{prob:2} and the corresponding Lax formula}
\label{sec:Cvx_Analysis_Prob2_TI}

As summarized in Table \ref{tab:comparisonCvxCond}, the temporally discretized Lax formula for Problem \ref{prob:2} ($\varphi_2$ in \eqref{eq:def_varphi2_disc} subject to \eqref{eq:def_varphi_const2_disc}) is non-convex for non-zero costs.
Thus, gradient-based methods do not, in general, provide a global optimal solution for Problem \ref{prob:2}.
Nevertheless, the Lax formula for the time-invariant Problem \ref{prob:2} provides a convex problem under a particular set of conditions, which will be presented in this subsection.

Consider a temporal discretization $\{t_0=0,...,T_K=T\}$, and the temporal discretized Lax formula $\varphi_2^\textnormal{TI}$ in \eqref{eq:def_varphi2_TI} subject to \eqref{eq:def_varphi2_TI_const}:
\begin{align}
        & \varphi_2^\textnormal{TI}(0,x) \simeq \min_{\mathrm{x}[\cdot],\beta[\cdot]}\sum_{k=0}^K H_2^{\textnormal{TI}*}(\mathrm{x}[k],\beta[k])\Delta_k+g(\mathrm{x}(K)), 
    \label{eq:def_varphi_TI2_disc}\\
    & \text{subject to } \begin{cases} \mathrm{x}[k+1]-\mathrm{x}[k] = -\Delta_k \beta[k],  \\ 
    \quad\quad\quad\quad\quad\quad\quad\quad\quad\quad\phantom{'''} k\in\{0,...,K-1\},\\
    \beta[k]\in \text{co}(\{0\}\cup B(\mathrm{x}[k])), \\
    \quad\quad\quad\quad\quad\quad\quad\quad\quad\quad\phantom{'''} k\in\{0,...,K-1\},\\
    \mathrm{x}[0]=x,\\
    c(\mathrm{x}[k]) \leq 0, \quad\quad\quad\quad\quad\phantom{'''} k\in \{0,...,K\}. \\  \end{cases}
    \label{eq:def_varphi_const_TI2_disc}
\end{align}
Note that $H_2^{\textnormal{TI}*}$ is defined in \eqref{eq:Ham_star2_TI_x}.

This problem is convex if the costs ($H_2^{\textnormal{TI*}}$ and $g$) in \eqref{eq:def_varphi_TI2_disc} and the constraints in \eqref{eq:def_varphi_const_TI2_disc} are convex in $(\mathrm{x}[\cdot],\beta[\cdot])$.
The last column of Table \ref{tab:comparisonCvxCond} shows the convexity result for the Lax formula in the time-invariant case.




\subsection{Comparison of convexity conditions for Problem \ref{prob:1}, \ref{prob:2}, and the corresponding Lax formulae}

Table \ref{tab:comparisonCvxCond} shows the summary of convexity conditions for Problems \ref{prob:1} and \ref{prob:2}, analyzed in Section \ref{sec:Cvx_Analysis_Prob1}, \ref{sec:Cvx_Analysis_Prob2}, and \ref{sec:Cvx_Analysis_Prob2_TI}.

In contrast to Problem \ref{prob:1}, the temporally discretized Lax formula for Problem \ref{prob:2} is generally non-convex unless the cost is zero.
This non-convexity is caused by the temporal discretization that converts $\min_{\tau\in[0,T]}$ in \eqref{eq:def_vartheta2} to $\min_{k'\in\{0,...,K\}}$ in \eqref{eq:def_vartheta2_disc2}.
However, this non-convexity issue is resolved in the temporally discretized Lax formula for the time-invariant case under the specified conditions in Table \ref{tab:comparisonCvxCond}.

\begin{remark}
\begin{enumerate}
    \item Suppose the stage cost $L$ and dynamics $f$ are decomposed into state and control dependent parts:
    \begin{align}
        & L(s,x,a) = L^x(s,x) + L^a(s,a),\\
        & f(s,x,a) = M(s)x + f^a(s,a).
    \end{align}
    $L^a$ and $f^a$ need to be convex and affine in the control, respectively, for convex temporally discretized Problem \ref{prob:1} ($\vartheta_1$ in \eqref{eq:def_vartheta1_disc} subject to \eqref{eq:def_vartheta_const1_disc}), but not for convex temporally discretized Lax formula $\varphi_1$ (\eqref{eq:def_varphi1_disc1} subject to \eqref{eq:def_varphi_const1_disc1}).
    \label{remark:benefit_prob1_Lax}
    
    \item The temporally discretized Problem \ref{prob:2} and Lax formula for Problem \ref{prob:2} are non-convex unless the cost is zero.
    
    \item Despite of the above remark, the temporally discretized Lax formula for the time-invariant Problem \ref{prob:2} ($\varphi_2^\textnormal{TI}$ in \eqref{eq:def_varphi_TI2_disc} subject to \eqref{eq:def_varphi_const_TI2_disc}) is convex if the stage cost and dynamics only depend on the control, and the terminal cost and state constraint are convex in the state.
\end{enumerate}
\end{remark}

    

    

\begin{table*}[ht]
    {
    \caption{Convexity Conditions for the Temporally Discretized Problems \ref{prob:1}, \ref{prob:2}, and the Corresponding Lax Formulae}
    \renewcommand{\arraystretch}{1.5}%
    \begin{tabular}{c||c|c||c|c|c}
        \hline\hline
          & Problem \ref{prob:1} & \begin{tabular}{c} Lax formula\vspace{-.2cm}\\for Problem \ref{prob:1}\end{tabular} & Problem \ref{prob:2} & \begin{tabular}{c} Lax formula\vspace{-.2cm}\\for Problem \ref{prob:2}\end{tabular} & \begin{tabular}{c} Lax formula\vspace{-.2cm}\\for time-invariant\vspace{-.2cm}\\Problem \ref{prob:2}\end{tabular} \\\hline\hline
         problem & $\vartheta_1$ \eqref{eq:def_vartheta1} s.t. \eqref{eq:def_vartheta_const1} & $\varphi_1$ \eqref{eq:def_varphi1} s.t. \eqref{eq:def_varphi_const1} & $\vartheta_2$ \eqref{eq:def_vartheta2} s.t. \eqref{eq:def_vartheta_const2} & $\varphi_2$ \eqref{eq:def_varphi2} s.t. \eqref{eq:def_varphi_const2} & $\varphi_2^\textnormal{TI}$ \eqref{eq:def_varphi2_TI} s.t. \eqref{eq:def_varphi2_TI_const}\\\hline
         \begin{tabular}{c} temporally\vspace{-.2cm}\\discretized\vspace{-.2cm}\\problem \end{tabular} & \eqref{eq:def_vartheta1_disc} s.t. \eqref{eq:def_vartheta_const1_disc} & \eqref{eq:def_varphi1_disc1} s.t. \eqref{eq:def_varphi_const1_disc1} & \eqref{eq:def_vartheta2_disc2} s.t. \eqref{eq:def_vartheta_const2_disc2} & \eqref{eq:def_varphi2_disc} s.t. \eqref{eq:def_varphi_const2_disc} & \eqref{eq:def_varphi_TI2_disc} s.t. \eqref{eq:def_varphi_const_TI2_disc} \\
         \hline\hline
         \begin{tabular}{c} stage cost \vspace{-.2cm}\\$L=L^x(s,x)$ \vspace{-.2cm}\\  \phantom{'''}$+L^a(s,a)$ \end{tabular} &
         \begin{tabular}{l} $L^x(s,\cdot)$ is \vspace{-.2cm}\\\phantom{'''''''''}convex in $x$  \vspace{-.2cm}\\$L^a(s,\cdot)$ is \vspace{-.2cm}\\\phantom{'''''''''}convex in $a$ \end{tabular}
          & $L^x(s,\cdot)$ is convex in $x$  & $L=0$ & $L=0$ & \begin{tabular}{c} $L^x=0$ \vspace{-.2cm}\\($L=L^a$) \end{tabular}\\\cline{1-3}\cline{6-6}\cline{3-5}
         \begin{tabular}{c} terminal cost \vspace{-.2cm}\\$g(s,x)$ \end{tabular} & convex in $x$& convex in $x$  & $g=0$ & $g=0$ & \begin{tabular}{c} $g=g(x)$ \vspace{-.2cm}\\is convex in $x$ \end{tabular}
          \\\cline{1-3}\cline{6-6}\cline{3-5}
         \begin{tabular}{c} dynamics \vspace{-.2cm}\\$f(s,x,a)$ \end{tabular}
          &\begin{tabular}{l} $f=M(s)x $ \vspace{-.2cm}\\\phantom{'''''''}$+N(s)a+C(s)$ \end{tabular} &\begin{tabular}{l} $f=M(s)x $ \vspace{-.2cm}\\\phantom{'''''''''}$+f^a(s,a)$ \end{tabular} &\begin{tabular}{l} $f=M(s)x $ \vspace{-.2cm}\\\phantom{'''''''}$+N(s)a+C(s)$ \end{tabular}  & \begin{tabular}{l} $f=M(s)x $ \vspace{-.2cm}\\\phantom{'''''''''}$+f^a(s,a)$ \end{tabular} & $f=f^a(a)$\\
         \cline{1-3}\cline{6-6}\cline{3-5}
         \begin{tabular}{c} state constraint \vspace{-.2cm}\\$c(s,x)$ \end{tabular} &  convex in $x$& convex in $x$ & $c=c(s)$ & $c=c(s)$ & convex in $x$\\
         \hline\hline
         \multicolumn{6}{p{502pt}}{For $s\in[0,T]$, $M$ (or $M(s)$) is a matrix in $\R^n\times\R^n$, $N$ (or $N(s)$) is a matrix in $\R^n\times\R^m$, $C(s)$ is a vector in $\R^n$, and $f^a(a)$ (or $f^a(s,a)$) is a nonlinear function in $a$, where $n$ and $m$ are the dimension of the state and the control.  }
    \end{tabular}
    }
    \label{tab:comparisonCvxCond}
\end{table*}

\section{Numerical Algorithm Using the Lax formula}
\label{sec:NumericalAlg}
\subsection{Numerical algorithm for Problem \ref{prob:1}}
\label{sec:numericalAlg_Prob12}


Remark \ref{remark:benefit_prob1_Lax} describes one of benefits of convexity of the Lax formula for Problem \ref{prob:1}.
For problems where Problem \ref{prob:1} is non-convex but the corresponding Lax formula is convex, we 1) solve the Lax formula for Problem \ref{prob:1} ($\varphi_1$ in \eqref{eq:def_varphi1} subject to \eqref{eq:def_varphi_const1}) by temporal discretization and gradient-based methods to get a global optimal solution $(\mathrm{x}_*,\beta_*)$ and 2) find the corresponding optimal state trajectory and control signal ($\mathrm{x}_*,\alpha_*$) to Problem \ref{prob:1} ($\vartheta_1$ in \eqref{eq:def_vartheta1} subject to \eqref{eq:def_vartheta_const1}).


In this section, we utilize a numerical method presented in \cite{donggunTAC_submitted} for a different class of problems, which also works for our problems.


Suppose an optimal solution $\mathrm{x}_*,\beta_*$ to the Lax formula for Problem \ref{prob:1} ($\varphi_1$) is given, in which the initial time $t$ is 0.
Consider a temporal discretization $\{t_0=0,...,t_K=T\}$ to numerically find approximate optimal stage trajectory $\mathrm{x}^\epsilon$ and control signal $\alpha^\epsilon$ for Problem \ref{prob:1} ($\vartheta_1$).

\cite{donggunTAC_submitted} showed that $H^*$ (the stage cost of the Lax formula for Problem \ref{prob:1}) is linearly decomposed by a finite number of stage costs.
This was derived from one of properties of the Legendre-Fenchel transformation: for $(s,x,b)\in[0,T]\times\R^n\times B(s,x)$,
\begin{align}
    H^*(s,x,b) = ((L^b)^*)^*(s,x,b),
    \label{eq:HamConj_Lb}
\end{align}
where $L^b (s,x,b)\coloneqq \min_{a\in A} L(s,x,a)$ subject to $f(s,x,a)=-b$, which is the corresponding stage cost for Problem \ref{prob:1} by converting the dynamics from $f(s,x,a)$ to $-b$.
This implies that, for all $k=0,...,K-1$, there exist $a_i^k\in A$ and $\gamma_i^k\in[0,1]$ ($\sum_i \gamma_i^k=1$) such that 
\begin{align}
\begin{split}
    \begin{bmatrix}H^*(t_k,\mathrm{x}_*(t_k),\beta_*(t_k))\\ -\beta_*(t_k)\end{bmatrix} = \sum_i \gamma_i^k \begin{bmatrix}L(t_k,\mathrm{x}_*(t_k),a_i^k)\\f(t_k,\mathrm{x}_*(t_k),a_i^k) \end{bmatrix}.
\end{split}
    \label{eq:cost_decomp}
\end{align} 
The detailed explanation for \eqref{eq:cost_decomp} can be found in \cite{donggunTAC_submitted}.
This implies that applying $\beta_*(t_k)$ for the time interval $[t_k,t_{k+1})$ in the Lax formula for Problem \ref{prob:1} ($\varphi_1$) can be approximated by applying a series of $a_i^k$ for $\gamma_i^k$-fraction of $[t_k,t_{k+1})$, which provides an approximate state trajectory and cost:
\begin{align}
    \alpha^\epsilon_*(s)=a_i^k, s\in [t_{k,i},t_{k,i+1}),
    \label{eq:approx_optCtrl_alpha}
\end{align}
where $t_{k,i}\coloneqq t_k+\sum_{j=1}^{i-1}\gamma_j^k \Delta_k$, $\Delta _k\coloneqq t_{k+1}-t_k$.
Then, $\mathrm{x}^\epsilon$ solving \eqref{eq:dynamics} for $s=0,\alpha=\alpha_*^\epsilon$ is an approximate optimal state trajectory for Problem \ref{prob:1}.

Theorem \ref{thm:OptCtrl_PostProcess} guarantees the quality of this approximation: as the size of the temporal discretization converges to 0, the approximation error of the state trajectory and the cost also converges to 0.
\begin{theorem}\textnormal{\textbf{(Guarantee of the quality of the approximation)}}
    Suppose Assumption \ref{assum:BigAssum} holds. For initial time $t=0$ and state $x\in\R^n$, assume that Problem \ref{prob:1} has a feasible solution. Denote an optimal control signal $\beta_*$ and an optimal state trajectory $\mathrm{x}_*$ to the Lax formula for Problem \ref{prob:1} ($\varphi_1$ in \eqref{eq:def_varphi1} subject to \eqref{eq:def_varphi_const1}).
    
    A control signal $\alpha_*^\epsilon$ in \eqref{eq:approx_optCtrl_alpha} on a temporal discretization $\{t_0=0,...,t_K=T\}$ and the corresponding state trajectory $\mathrm{x}_*^\epsilon$ solving \eqref{eq:dynamics} for $t=0,\alpha=\alpha_*^\epsilon$ satisfy the following statement:
    let $\delta = \max_{k}\Delta_k$, then
    \begin{align}
        &  \lim_{\delta \rightarrow 0} \| \mathrm{x}_* - \mathrm{x}_*^\epsilon \|_{L^\infty (0,T;\R^n)} =0, \label{eq:thm_OptCtrlProcess_eq1}\\
        &  \lim_{\delta \rightarrow 0}
        \| c(\cdot,\mathrm{x}_*^\epsilon(\cdot)) - c(\cdot,\mathrm{x}_*(\cdot)) \|_{L^\infty (0,T)} =0,
        \label{eq:thm_OptCtrlProcess_eq2}\\
        \begin{split}
            &\lim_{\delta \rightarrow 0} \max_{\tau\in[0,T]}\int_0^\tau L(s,\mathrm{x}_*^\epsilon(s),\alpha_*^\epsilon(s)) ds+g(\tau,\mathrm{x}_*^\epsilon(\tau)) \\
            &\quad\quad\quad\quad  = \vartheta_1(0,x),    
        \end{split}
        \label{eq:thm_OptCtrlProcess_eq3}
    \end{align}
    where $\vartheta_1$ is defined in \eqref{eq:def_vartheta1} subject to \eqref{eq:def_vartheta_const1}.
    \label{thm:OptCtrl_PostProcess}
\end{theorem}
\textbf{Proof.} See Appendix \ref{appen:thm_OptCtrl_PostProcess}.

Algorithm \ref{alg:Opt_NewFormulation1} summarizes how to numerically solve Problem \ref{prob:1} ($\vartheta_1$) using the Lax formula $\varphi_1$ with the temporal discretization $\{t_0=0,...,t_K=T\}$.

\begin{algorithm}[t]
\caption{Computing optimal state trajectory ($\mathrm{x}_*^\epsilon$) and control signal ($\alpha_*^\epsilon$) for Problem \ref{prob:1} using the Lax formula}
\begin{algorithmic}[1]
\State \textbf{Input:} {initial time $t=0$, initial state $x$.}
\State \textbf{Output:} {optimal state trajectory ($\mathrm{x}$), control signal ($\alpha$), and terminal time ($\tau$).}
\State Generate a temporal discretization: $\{t_0=t,...,t_K=T\}$.
\State Solve $\varphi_1$ in \eqref{eq:def_varphi1_disc1} subject to \eqref{eq:def_varphi_const1_disc1} for $\mathrm{x}_* [\cdot],\beta_*[\cdot]$.
\State Solve $(a_i^k,\gamma_i^k)$ in \eqref{eq:cost_decomp}.
\State Additionally discretize each temporal interval $[t_k,t_{k+1})$ into multiple sub intervals: $[t_{k,i},t_{k,i+1})$, and design $\alpha_*^\epsilon$ as in \eqref{eq:approx_optCtrl_alpha}.
\State Compute $\mathrm{x}_*^\epsilon$ by solving the ODE \eqref{eq:dynamics} for $\alpha_*^\epsilon$ and the initial state $x$.
\State Compute $\tau_*=t_{k'}$, where $k'$ maximizes \eqref{eq:def_vartheta1_disc}.
\end{algorithmic}
\label{alg:Opt_NewFormulation1}
\end{algorithm}

With the analogous derivation, one can prove that the approximate optimal control design in \eqref{eq:approx_optCtrl_alpha} also works for Problem \ref{prob:2} and the corresponding Lax formula. 
Nevertheless, these two problems are generally non-convex, thus, gradient-based methods do not provide a globally optimal solution $(\mathrm{x}_*,\beta_*)$ for the temporally discretized Lax formula for Problem \ref{prob:2}.
The following subsection presents a numerical method for the time-invariant Problem \ref{prob:2}, in which the corresponding Lax formula is convex under the conditions in the last column of Table \ref{tab:comparisonCvxCond}.

\subsection{Numerical algorithm for the time-invariant Problem \ref{prob:2}}
\label{sec:NumericalAlg_TI_Prob2}

We present a numerical method for the time-invariant Problem \ref{prob:2} based on one of the authors' prior work \cite{lee2020hopf}, which provides a numerical method for the reach-avoid problem, which is a particular example of Problem \ref{prob:2} where the stage cost is zero $L=0$.
We generalize the method \cite{lee2020hopf} for non-zero stage cost even though the convexity condition for the time-invariant Problem \ref{prob:2} in Table \ref{tab:comparisonCvxCond} requires zero stage cost.

Suppose an optimal solution $\mathrm{x}_*,\beta_*$ to the Lax formula for the time-invariant Problem \ref{prob:2} is given for the zero initial time, $t=0$, and consider the temporal discretization $\{t_0=0,...,t_K=T\}$.

In \cite{donggunTAC_part1}, the time-invariant Problem \ref{prob:2} is converted to a fixed-horizon optimal control problem, in which the stage cost is $L(x,a)a_d$ and $a_d\in[0,1]$ is a freezing control. 
In the fixed-horizon optimal control problem, the dynamics $f(x,a)a_d$ is converted to $-b$, and the stage cost $L$ is also converted to $L^{b,\textnormal{TI}}(x,b)\coloneqq\min_{a\in A,a_d \in[0,1]}L(x,a)a_d$ subject to $f(x,a)a_d=b$.
Similar to \eqref{eq:HamConj_Lb}, we have, for $(x,b)\in \R^n\times B(x)$,
\begin{align}
    H_2^{\textnormal{TI}*}(x,b) = (L^{b,\textnormal{TI}})^{**}(x,b).
\end{align}
This implies that the epigraph of $H_2^{\textnormal{TI}*}(x,\cdot)$ \eqref{eq:def_ConvConj_Ham} is the convex hull of the epigraph of $L^{b,\textnormal{TI}}(x,\cdot)$, and the domain of $H_2^{\textnormal{TI}*}(x,\cdot)$ is the convex hull of the union of the zero vector and $B(x)$ in \eqref{eq:ctrlbound_B}.
Thus, for all $k=0,...,K-1$, there exist $a_i^k\in A$ and $\gamma_i^k\in[0,1]$ ($\sum_i \gamma_i^k \leq 1$) such that 
\begin{align}
\begin{split}
    \begin{bmatrix}H_2^{\textnormal{TI}*}(\mathrm{x}_*(t_k),\beta_*(t_k))\\ -\beta_*(t_k)\end{bmatrix} = \sum_i \gamma_i^k \begin{bmatrix}L(\mathrm{x}_*(t_k),a_i^k)\\f(\mathrm{x}_*(t_k),a_i^k) \end{bmatrix}.
\end{split}
    \label{eq:cost_decomp2}
\end{align} 
 
We approximate the state trajectory and control signal $\mathrm{x}_*^\epsilon, \beta_*^\epsilon$ in two steps.
In step 1, we first define approximate $\mathrm{x}_*^1,\beta_*^1$:
\begin{align}
    \beta_*^1(s)=\begin{cases} -f(x_*(t_k),a_i^k),&s\in[t_{k,i},t_{k,i+1}),\\ 0,& s\in[t_k+\sum_{i}\gamma_i^k\Delta_k,t_{k+1}),\end{cases}
    \label{eq:ctrl_beta_star_1}
\end{align}
where $t_{k,i}\coloneqq t_k+\sum_{j=1}^{i-1}\gamma_j^k \Delta_k$, and $\mathrm{x}_*^1$ solves
\begin{align}
    \dot{\mathrm{x}}_*^1(s)=-\beta_*^1(s),\quad s\in[0,T],\quad \mathrm{x}_*^1(0)=x.
\end{align}
In step 2, the psuedo time operation is introduced.
\begin{align}
    \sigma_{\beta_*^1}(s)\coloneqq \int_0^s \mathds{1}\{ \beta_*^1(\tau)\neq 0 \} d\tau, 
    \label{eq:psuedoTime}
\end{align}
where $\mathds{1}$ outputs 1 if the input condition is true, or 0 if not. 
The corresponding inverse operation is defined: for $s\in[0,T]$,
\begin{align}
    \sigma^{-1}_{\beta_*^1}(s) \coloneqq \min \tau \text{ subject to } \sigma_{\beta_*^1}(\tau ) = s.
\end{align}
Using this operation, define state trajectory and control signal: 
\begin{align}
    \alpha_*^\epsilon (s) = \begin{cases} a_i^k, &s\in[\sigma_{\beta_*^1}(t_i^k),\sigma_{\beta_*^1}(t_{i+1}^k)),\\ \text{any control}, & s\in[\sigma_{\beta_*^1}(T),T], \end{cases}
    \label{eq:approx_optCtrl_alpha2}
\end{align}
and $\mathrm{x}_*^\epsilon$ solves \eqref{eq:dynamics} for $\alpha_*^\epsilon$ with the initial state $x$.
Theorem \ref{thm:OptCtrl_PostProcess2} shows $(\mathrm{x}_*^\epsilon,\alpha_*^\epsilon)$ are approximate optimal state trajectory and control signal to the time-invariant Problem \ref{prob:2}.

\begin{theorem}\textnormal{\textbf{(Guarantee of the quality of the approximation for the time-invariant Problem \ref{prob:2})}}
    Suppose Assumption \ref{assum:BigAssum} holds,
    and assume that Problem \ref{prob:2} for the time-invariant case has a feasible solution. For initial time $t=0$ and state $x\in\R^n$, denote an optimal control signal $\beta_*$ and an optimal state trajectory $\mathrm{x}_*$ to the Lax formula for the time-invariant Problem \ref{prob:2} ($\varphi_2^\textnormal{TI}$ in \eqref{eq:def_varphi2_TI} subject to \eqref{eq:def_varphi2_TI_const}).
    
    A control signal $\alpha_*^\epsilon$ in \eqref{eq:approx_optCtrl_alpha2} on a temporal discretization $\{t_0=0,...,t_K=T\}$ and the corresponding state trajectory $\mathrm{x}_*^\epsilon$ solving \eqref{eq:dynamics} for $t=0,\alpha=\alpha_*^\epsilon$ satisfy the following statement:
    let $\delta = \max_{k}\Delta_k$, then
    \begin{align}
        &\lim_{\delta\rightarrow0}\lVert \mathrm{x}_* (\cdot)- \mathrm{x}_*^\epsilon(\sigma(\cdot;\beta_*^1))  \rVert_{L^\infty (0,T;\R^n)} =0,\label{eq:thm_OptCtrlProcess2_eq1}\\
        &\lim_{\delta\rightarrow0}\lVert c(\cdot,\mathrm{x}_* (\cdot))- c(\cdot, \mathrm{x}_*^\epsilon(\sigma(\cdot;\beta_*^1)))  \rVert_{L^\infty (0,T)} =0,
        \label{eq:thm_OptCtrlProcess2_eq2}\\
        \begin{split}
        &\lim_{\delta\rightarrow0}  \int_0^{\sigma(T;\beta_*^1)} L(\mathrm{x}_*^\epsilon(s),\alpha_*^\epsilon(s)) ds  +g\left(\mathrm{x}_*^\epsilon(\sigma(T;\beta_*^1))\right)\\ &\quad\quad\quad\quad =\vartheta_2(0,x),
    \end{split}
        \label{eq:thm_OptCtrlProcess2_eq3}
    \end{align}
    where $\vartheta_2$, $\sigma$, and $\beta_*^1$ are defined in \eqref{eq:def_vartheta2} subject to \eqref{eq:def_vartheta_const2}, \eqref{eq:psuedoTime}, and \eqref{eq:ctrl_beta_star_1}.
    \label{thm:OptCtrl_PostProcess2}
\end{theorem}
\textbf{Proof.} See Appendix \ref{appen:thm_OptCtrl_PostProcess2}.

Algorithm \ref{alg:Opt_NewFormulation2} summarizes how to numerically solve the time-invariant Problem \ref{prob:2} ($\vartheta_2$) using the Lax formula $\varphi_2^\textnormal{TI}$ with the temporal discretization $\{t_0=0,...,t_K=T\}$.

\begin{algorithm}[t]
\caption{Computing optimal state trajectory ($\mathrm{x}_*^\epsilon$) and control signal ($\alpha_*^\epsilon$) for the time-invariant Problem \ref{prob:2} using the Lax formula 
}
\begin{algorithmic}[1]
\State \textbf{Input:} {initial time $t=0$, initial state $x$}.
\State \textbf{Output:} {optimal state trajectory ($\mathrm{x}$), control signal ($\alpha$), and terminal time ($\tau$)}.
\State Generate a temporal discretization: $\{t_0=t,...,t_K=T\}$.
\State Solve $\varphi_2^\textnormal{TI}$ in \eqref{eq:def_varphi_TI2_disc} subject to \eqref{eq:def_varphi_const_TI2_disc} for $\mathrm{x}_* [\cdot],\beta_*[\cdot]$ .
\State Solve $(a_i^k,\gamma_i^k)$ in \eqref{eq:cost_decomp2}.
\State Additionally discretize each temporal interval $[t_k,t_{k+1})$ into multiple sub intervals: $[t_{k,i},t_{k,i+1})$, and design $\beta*^1$ as in \eqref{eq:ctrl_beta_star_1}.
\State Design $\alpha_*^\epsilon$ \eqref{eq:approx_optCtrl_alpha2}
\State Compute $\mathrm{x}_*^\epsilon$ by solving the ODE \eqref{eq:dynamics} for $\alpha_*^\epsilon$ and the initial state $x$.
\State Compute $\tau_*=t_{k'}$, where $k'$ minimizes \eqref{eq:def_vartheta2_disc2}.
\end{algorithmic}
\label{alg:Opt_NewFormulation2}
\end{algorithm}

\section{Numerical Examples} 
\label{sec:example}

This section provides two numerical examples to demonstrate the Lax formulae for Problems \ref{prob:1} and \ref{prob:2}.

\subsection{Problem \ref{prob:1}: Robust formation control}

We control a multi-robot inspection system.
Each robot moves in two-dimensional space and is equipped with a sensor system.
The robots need to maintain a certain formation during the sensing time, so that each sensor image and the gathered image map are of good quality.
In the presence of disturbance, the formation might not be maintained, and we would like to robustly control the multi-robot system so that the formation violation is minimized. 
Suppose the required sensing time is 2 s, and each robot is a four dimensional system, which follows 
\begin{align}
\begin{split}
    & \dot{\mathrm{x}}_1^r(s) = \mathrm{x}_2^r (s),\quad \dot{\mathrm{x}}_2^r(s) = \alpha_1^r (s) \cos \alpha_2^r (s)+d(s),\\
    & \dot{\mathrm{x}}_3^r(s) = \mathrm{x}_4^r (s),\quad \dot{\mathrm{x}}_4^r(s) = \alpha_1^r (s) \sin \alpha_2^r (s),
\end{split}
\label{eq:exaple1_dyn}
\end{align}
where $\mathrm{x}^r(s)=(\mathrm{x}^r_1(s),\mathrm{x}^r_2(s),\mathrm{x}^r_3(s),\mathrm{x}^r_4(s))\in \R^4$ for all $s\in[0,2]$, $\mathrm{x}=(\mathrm{x}^1,...,\mathrm{x}^R)\in\R^{4R}$, $R=10$ is the number of the robots, $r=1,...,R$, and $d$ is the horizontal disturbance: $d(s)=0.5(1+\cos \pi s)$.
$\mathrm{x}_1$ and $\mathrm{x}_2$ ($\mathrm{x}_3$ and $\mathrm{x}_4$) are horizontal (vertical) position and velocity, $\alpha_1^r$ is the magnitude of acceleration, and $\alpha_2^r$ is the angle of the acceleration of robot $r$.

Consider the following problem:
\begin{align}
\begin{split}
    &\quad\quad  \inf_{\alpha}\max_{\tau\in[0,2]} \lVert (\mathrm{x}^1_1(\tau),\mathrm{x}^1_3(\tau))-x_g^1\rVert_2 \\
    &\quad\quad +\sum_{r= 2}^{R} \lVert (\mathrm{x}^r_1(\tau),\mathrm{x}^i_3(\tau)) - (\mathrm{x}^r_1(\tau),\mathrm{x}^1_3(\tau)) -o^r \rVert_2
    \label{eq:example1_prob_cost}
\end{split}\\
&\text{subject to} \begin{cases} \eqref{eq:exaple1_dyn},\text{ } \alpha^r(s)\in[-1,1]\times[-\frac{\pi}{6},\frac{\pi}{6}],\text{ } s\in[0,2],& \\\mathrm{x}(0)=x,\\  \mathrm{x}_1^r(s)\leq 5.2, \mathrm{x}_2^r(s)\geq 0,\quad r=1,...,R, s\in[0,2]. \end{cases}
\label{eq:example1_prob_const}
\end{align}
The stage cost $L$ is zero, $x_g^1$ is $(1,1)$, $o^r=(0.4rR/(R-1),0)(\in\R^2)$ is the offset for the formation of robot $i$ with respect to robot $1$, and the initial state of the robots is the positions whose corresponding cost is zero with randomly chosen velocities between -0.5 and 0.5: $(\mathrm{x}^1_1(0),\mathrm{x}^1_3(0))=x_g^1$ and $(\mathrm{x}^r_1(\tau),\mathrm{x}^r_3(\tau)) = (\mathrm{x}^1_1(\tau),\mathrm{x}^1_3(\tau)) +o^r$. 

For the given problem, the Hamiltonian $H$ in \eqref{eq:def_Ham} becomes
\begin{align}
    &H(s,x,p)=\bar{H}(s,x,z,p,-1) \notag\\
    =& \max_{\substack{a^r\in [-1,1]\times[-\frac{\pi}{6},\frac{\pi}{6}], \\ r=1\dots R}} - \sum_{r} \begin{tabular}{l} $\big( p^r_1 x^r_2 + p^r_2  (a_1^r \cos a_2^r + d(s))$\\
    $+ p^r_3 x_4^r  + p^r_4 a_1^r\sin a_2^r \big)$\end{tabular} \notag\\
    =& \sum_{r} \begin{tabular}{l}$-p_1^r x_2^r -p_2^r d(s) - p_3^r x_4^r$\\ $+ \max\{ \|(p_2^r,p_4^r)\|_2, |p_2^r|\frac{\sqrt{3}}{2} + |p_4^r|\frac{1}{2} \}$.\end{tabular}
\end{align}
$H$ is convex in $p$, and any supporting hyperplane for $H$ can be written as $b\cdot p=0$ for some normal vector $b$. 
Since the supporting hyperplane $b\cdot p=0$ passes through the origin for any $b$,
\begin{align}
    H^*(s,x,b)=0
\end{align}
for $b\in\text{Dom}(H^*(s,x,\cdot))=\text{co}(B(s,x))$.
See Lemma \ref{lemma:beta_const} for the above property for the domain of $H^*$. 
In general, if the stage cost $L$ is zero, $H^*$ becomes zero, which has been investigated in \cite{donggunTAC_submitted}.

By the Lax formula for Problem \ref{prob:1} in Theorem \ref{thm:HJeq_prob1}, \eqref{eq:example1_prob_cost} subject to \eqref{eq:example1_prob_const} is equivalent to 
\begin{align}
\begin{split}
    &\quad\quad  \inf_{\beta}\max_{\tau\in[0,2]} \lVert (\mathrm{x}^1_1(\tau),\mathrm{x}^1_3(\tau))-x_g^1\rVert_2 \\
    &\quad\quad +\sum_{i= 2}^{r} \lVert (\mathrm{x}^r_1(\tau),\mathrm{x}^r_3(\tau)) - (\mathrm{x}^1_1(\tau),\mathrm{x}^1_3(\tau)) - o^r \rVert_2
    \label{eq:example1_prob_Lax_cost}
\end{split}\\
&\text{subject to} \begin{cases} \mathrm{x}^r(s)=-\beta^r(s),\\
\beta^r_1(s)=-\mathrm{x}^r_2(s),\beta^r_3(s)=-\mathrm{x}^r_4(s),\\
 \lVert(\beta^r_2(s)+d(s),\beta^r_4(s))\rVert_2\leq 1, |\beta^r_4(s)|\leq \frac{1}{2},& \\
 \mathrm{x}(0)=x,\\\mathrm{x}_1^r(s)\leq 5.2, \mathrm{x}_2^r(s)\geq 0, \quad r=1,...,R, s\in[0,2], \end{cases}
\label{eq:example1_prob_Lax_const}
\end{align}
where $\beta(s)=(\beta^1,...,\beta^R)(s)$ and $\beta^r(s)=(\beta_1^r,\beta_2^r,\beta_3^r,\beta_4^r)(s)$.
Note that $B(s,x)=B^1(s,x^1)\times\dots\times B^R(s,x^R)$, 
\begin{align}
    B^r(s,x^r)=\left\{(-x_2^r,b_2,-x_4^r,b_4)\bigg|\begin{tabular}{l}$\lVert(b_2+d(s),b_4)\rVert_2\leq1,$\\ $|b_4|\leq\sqrt{3}|b_2 + d(s)| $\end{tabular} \right\},
\end{align}
and $\text{co}(B(s,x))$ can be found in \eqref{eq:example1_prob_Lax_const}.

The given problem \eqref{eq:example1_prob_cost} subject to \eqref{eq:example1_prob_const} is non-convex, but the Lax formula \eqref{eq:example1_prob_Lax_cost} subject to \eqref{eq:example1_prob_Lax_const} is convex since Condition \ref{cond:cond1_varphi1} is satisfied, which allows the proposed Lax formula to provide an optimal solution.

For numerical computation of the Lax formula, we discretize the temporal space to $\{t_0=0 \dots t_K = 2\}$ with $\Delta_k=0.1$ (21 time steps). The computation time to solve the Lax formula (\eqref{eq:example1_prob_Lax_cost} subject to \eqref{eq:example1_prob_Lax_const}) is 103.7 s, in which the optimal control signal $\beta_*$ and state trajectory $\mathrm{x}_*$ are computed.
This system is 40 dimension with 10 robots, for which it is intractable to utilize grid-based methods (such as the level-set method \cite{Osher02}) to solve the HJ equations \eqref{eq:HJeq1}.


We follow Algorithm \ref{alg:Opt_NewFormulation1} to compute an optimal control signal ($\alpha_*^\epsilon$ in \eqref{eq:approx_optCtrl_alpha}) for the given problem.
Applying line 5 to 7 in Algorithm \ref{alg:Opt_NewFormulation1} for each robot also works. 
In other words, we first find $a_i^{k,r}$ and $\gamma_i^{k,r}$ for robot $r=1,...,R$ as in \eqref{eq:cost_decomp}. 
For $\beta_*^r(t_k)=(\beta_{*,1}^r,\beta_{*,2}^r,\beta_{*,3}^r,\beta_{*,4}^r)(t_k)\in B^r(t_k,\mathrm{x}_*^r(t_k))$,
\begin{align}
    a_1^{k,r}=\begin{small}\begin{pmatrix}\text{sign}\left(-\beta_{*,2}^r(t_k)-d(t_k)\right)\left\|\begin{pmatrix} \beta_{*,2}^r(t_k)+d(t_k) \\\beta_{*,4}^r(t_k) \end{pmatrix}\right\|_2 \\\arctan\left(\beta_{*,4}^r(t_k)/\left(\beta_{*,2}^r(t_k)+d(t_k)\right) \right) \end{pmatrix}\end{small},
\end{align}
$\gamma_1^{k,r} = 1$, and for $\beta_*^r(t_k)$ in $\text{co}(B^r(t_k,\mathrm{x}_*^r(t_k)))$ but not in $B^r(t_k,\mathrm{x}_*^r(t_k))$, 
\begin{align}
    & a_1^{k,r} = \left( -2\beta_{*,4}^r(t_k), \pi/6 \right), a_2^{k,r} = -a_1^{k,r} ,\\
    & \gamma_1^{k,r} =\frac{|-\beta_*(t_k,2)-d(t_k)-\sqrt{3}\beta_{*,4}(t_k)|}{|2\sqrt{3}\beta_{*,4}(t_k)|},\\
    & \gamma_2^{k,r} =\frac{| \beta_{*,2}(t_k)+d(t_k)-\sqrt{3}\beta_{*,4}(t_k)|}{|2\sqrt{3}\beta_{*,4}(t_k)|}.
\end{align}
Based on the computed $a_1^{k,r}$ and $\gamma_i^{k,r}$ for robot $r$, each approximate optimal control signal $\alpha_*^{r,\epsilon}$ is designed as in \eqref{eq:approx_optCtrl_alpha}, and the corresponding state trajectory $\mathrm{x}_*^{\epsilon}$ can be also computed by solving the ODE \eqref{eq:dynamics} for $\alpha_*^{\epsilon}$, where $\alpha_*^{\epsilon}=(\alpha_*^{1,\epsilon},...,\alpha_*^{R,\epsilon})$.

\begin{figure}[t!]
\centering
\begin{tabular}{c}
    \includegraphics[trim = 0mm 0mm 0mm 0mm, clip, width=0.48\textwidth]{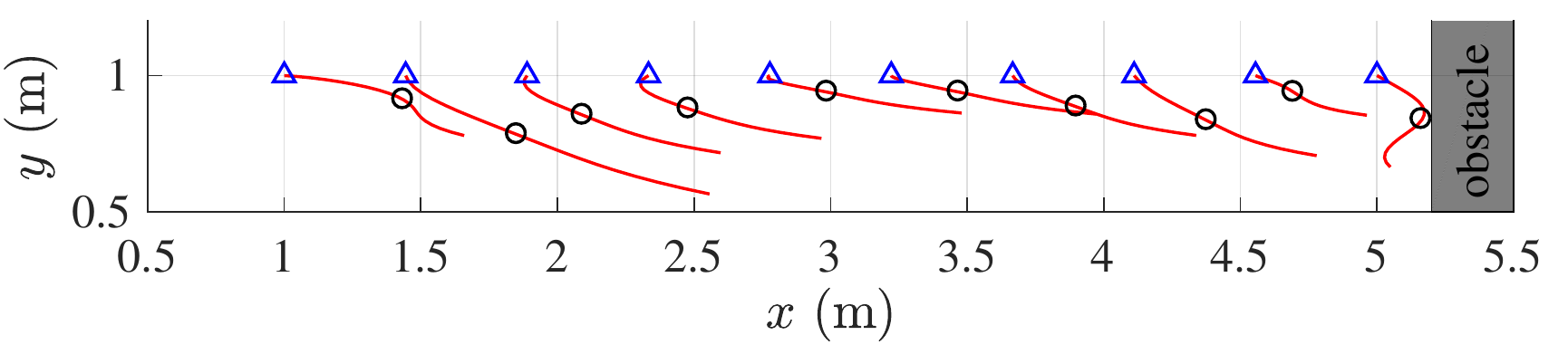}   \\
    \includegraphics[trim = 0mm 0mm 0mm 0mm, clip, width=0.40\textwidth]{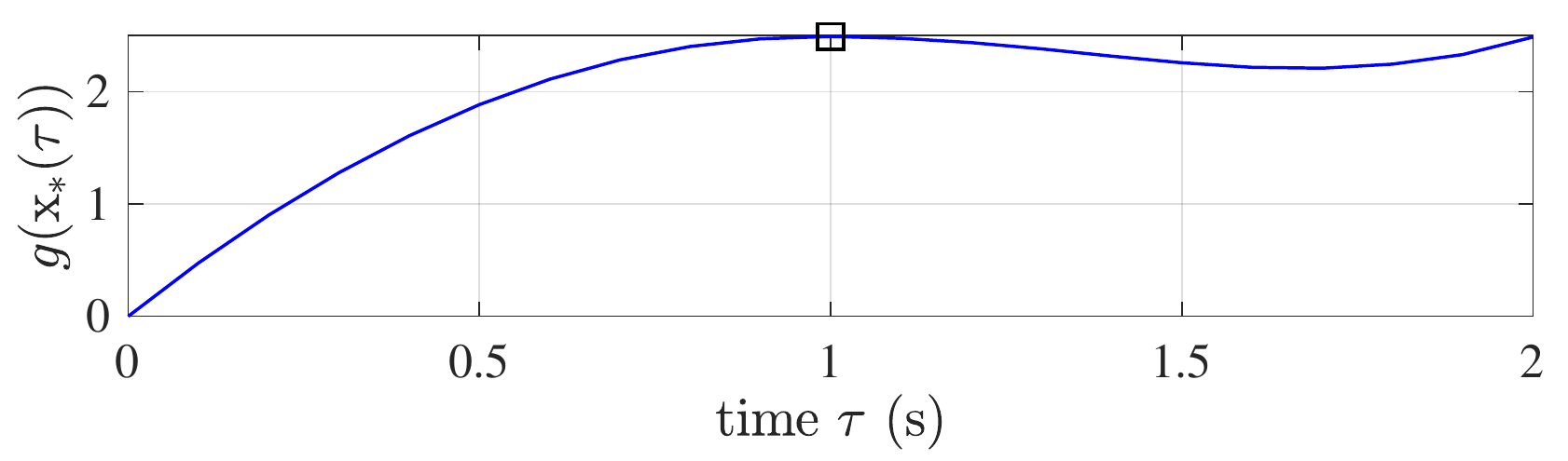} 
\end{tabular}
\caption{ (a) The red lines are the optimal state trajectories of the 10 robots that robustly minimize the worst formation violation in the given time horizon. The blue triangles are the initial position of the robots, and the black circles are the position of the robots at which the cost $g(\mathrm{x}_*(\tau))$ is maximized. (b) Under the optimal control, the worst cost is attained at 1 s.
}
\label{fig:example1_traj}
\end{figure}

Figure \ref{fig:example1_traj} (a) shows an optimal position of the robots, and Figure \ref{fig:example1_traj} (b) shows the cost over the time. 
The cost is maximized at 1 s, at which the positions of the robots are marked as the black circles in Figure \ref{fig:example1_traj} (a). 
The blue triangles are the initial positions of the robots.
Even under the disturbance, all robots avoid the obstacles.

\subsection{Lax formula for the time-invariant Problem \ref{prob:2}}


In this example, we also deal with a multi-robot system with each robot moving in two-dimensional space, and the goal is to find an optimal control signal and time that achieve the best formation of the robot system over the given time horizon.
Suppose the given time-horizon is 2 s, and each robot is a two dimensional system:
\begin{align}
    \dot{\mathrm{x}}_1^r(s) = \alpha_1^r (s) + 2, \quad \dot{\mathrm{x}}_2^r(s) = \alpha_2^r (s),
    \label{eq:ex2_dyn}
\end{align}
where $\mathrm{x}=(\mathrm{x}^1,...,\mathrm{x}^R)$, $R=10$ is the number of the robots, $\mathrm{x}^r(s)=(\mathrm{x}_1^r(s),\mathrm{x}_2^r(s))\in\R^2$ for all $s\in[0,2]$, $r=1,..,R$, $\alpha^r(s)=(\alpha_1^r,\alpha_2^r)(s)\in\R^2$, and $\|\alpha^r(s)\|_\infty \leq 1$.

Consider the following problem:
\begin{align}
    &\inf_{\alpha}\min_{\tau\in[0,2]}  \lVert \mathrm{x}^1(\tau)-x_g^1\rVert_2 + \sum_{r = 2}^{R} \lVert \mathrm{x}^r(\tau) - \mathrm{x}^1(\tau) - o^r\rVert_2
    \label{eq:ex2_prob_cost}
    \\
    &\text{subject to}\begin{cases}\eqref{eq:ex2_dyn}, \lVert \alpha^r(s)\rVert_\infty \leq 1, s\in[0,2],\\ \mathrm{x}(0)=x,\\
    \mathrm{x}_1^r(s) \leq 5, \mathrm{x}_2^r(s)\geq 0, r=1\dots R,s\in[0,\tau], \end{cases}
    \label{eq:ex2_prob_const}
\end{align}
where $x_g^1 = (1,1)$, $o^r=(0.4rR/(R-1),0)\in\R^2$ is the offset for the formation of robot $r$ with respect to robot $1$, the initial state of the first robot is randomly sampled near the goal state $x_g^1$, and for the other robots, the initial states are randomly sampled near $x_g^1+o^r$.

For this problem, the Hamiltonian in \eqref{eq:def_Hambar} is
\begin{align}
    \bar{H}(x,z,p,q) &= \max_{\substack{\|a^r\|_\infty\leq 1, \\ r=1\dots R}} -\sum_r \left(p^r_1 (a^r_1 + 2) + p^r_2 a^r_2\right)  \notag\\
    & = -2\sum_i p_1^r + \|p\|_1,
    \label{eq:ex2_eq2}
\end{align}
where $p=(p^1,...,p^R)\in\R^{2R}$, and $p^r=(p_1^r,p_2^r)\in\R^2$. 
By \eqref{eq:def_Ham} and \eqref{eq:Ham_Ham2TI}, 
\begin{align}
    H_2^\textnormal{TI}(x,p) = \max\{0,-2\sum_r p_1^r + \|p\|_1\}.
    \label{eq:ex2_eq3}
\end{align}
Since $H_2^\textnormal{TI}$ is a pointwise maximum of two convex functions in $p$, $H_2^\textnormal{TI}$ is convex in $p$. 
Also, for all $b\in\R^{2R}$, the supporting hyperplane of $H_2^\textnormal{TI}$ in $p$-space with respect to the normal vector $b$ crosses the origin. Thus,
\begin{align}
    H_2^{\textnormal{TI}*}(x,b) = 0
    \label{eq:ex2_eq4}
\end{align}
for $b=(b^1,...,b^R)=(b^1_1,b^1_2,...,b^R_1,b^R_2)\in\text{Dom}(H_2^{\textnormal{TI}*}(x,\cdot))=\text{co}(\{0\}\cup B(x))$, where $B(x)=\{[-3,-1]\times[-1,1]\}^R$ and
\begin{align}
     \text{co}&(\{0\}\cup B(x)) = \big\{  b ~|~\forall r_1 ,r_2\in\{1,...,R\},\notag\\
     & -3\leq b^{r_1}_1, |b_2^{r_1}|\leq 1, b^{r_1}_1-3b^{r_2}_1 \geq 0, b^{r_1}_1 -\frac{1}{3}b^{r_2}_1\leq 0, \notag\\
    & b^{r_1}_2 - b^{r_2}_1 \geq 0, b^{r_1}_2 + b^{r_2}_1 \leq 0 \big\}.
\label{eq:ex2_CtrlSet}
\end{align}

By the Lax formula for the time-invariant Problem \ref{prob:2} in Theorem \ref{thm:Laxformula_Prob2_TI} and \eqref{eq:ex2_eq4}, the given problem is equivalent to 
\begin{align}
    &\quad\quad  \inf_{\beta} \lVert \mathrm{x}^1 (2)-x_g^1\rVert_2 + \sum_{r=2}^R \lVert \mathrm{x}^r(2) - \mathrm{x}^1(2) - o^r \rVert_2
    \label{eq:ex2_Lax_cost}
\\
&\text{subject to }\begin{cases} \dot{\mathrm{x}}^r(s)=-\beta^r(s),\\
\beta(s)\in\text{co}\left(\{0\}\cup B(\mathrm{x}_*(s))\right) \text{ in \eqref{eq:ex2_CtrlSet}}, \\
\mathrm{x}(0)=x\\
\mathrm{x}_1^r(s) \leq 5, \mathrm{x}_2^r(s)\geq 0, r=1\dots R,s\in[0,2].
\end{cases} 
\label{eq:ex2_Lax_const}
\end{align}
The temporally discretized Lax formula for the time-invariant Problem 2 (\eqref{eq:ex2_Lax_cost} subject to \eqref{eq:ex2_Lax_const}) is convex although the temporally discretized given problem (\eqref{eq:example1_prob_cost} subject to \eqref{eq:example1_prob_const}) is non-convex.
Thus, the Lax formula (\eqref{eq:ex2_Lax_cost} subject to \eqref{eq:ex2_Lax_const}) provides an optimal solution by utilizing gradient-based methods. 
This formula is numerically solved by the interior-point method, in which numerical optimal control signal $\beta_*$ and state trajectory $\mathrm{x}_*$ are specified.
The computational time is 82.6 s.
Since the dimension of the system is 20, it is intractable to utilize any grid-based method to solve the HJ equation \eqref{eq:HJeq2_TI}.

Using $\beta_*$ and $\mathrm{x}_*$, we follow the steps in Algorithm \ref{alg:Opt_NewFormulation2}.
Denote $\beta_*(t_k)=(\beta_*^1,...,\beta_*^R)(t_k)$ and $\beta_*^r=(\beta_{*,1}^r,\beta_{*,2}^r)$.
First, we find $a_i^{k}$ and $\gamma_i^{k}$ as in \eqref{eq:ctrl_beta_star_1}.
For $\beta_*(t_k)\in B(\mathrm{x}_*(t_k))=\{[-3,-1]\times[-1,1]\}^R$, 
$a_1^k$ is $(a_1^{k,1},...,a_1^{k,r})\in\R^{2R}$, where
\begin{align}
    a_1^{k,r} = (-\beta_{*,1}^1(t_k)-2 , -\beta_{*,2}^1(t_k)), \text{ and }\gamma_1^k=1.
\end{align}
For $\beta_*(t_k)$ in $\text{co}(\{0\}\cup B(\mathrm{x}_*(t_k))\}$ but not in $B(\mathrm{x}_*(t_k))$, $a_1^k$ is  $(a_1^{k,1},...,a_1^{k,r})\in\R^{2R}$,
\begin{align}
    & a_1^{k,r} = (-\beta_{*,1}^r(t_k)/\gamma_1^k -2 , -\beta_{*,2}^r(t_k)/\gamma_1^k ),\\
    & \gamma_1^k = -\max_{r=1,...,R} \beta_{*,1}^r \in[0,1]\subset\R.
\end{align}
Using the above $a_i^{k,r}$ and $\gamma_i^{k,r}$, an approximate optimal control signal $\alpha_*^\epsilon$ is designed as in \eqref{eq:approx_optCtrl_alpha2}, and the corresponding state control trajectory $\mathrm{x}_*^\epsilon$ is computed by solving the dynamics ODE \eqref{eq:dynamics}.

\begin{figure}[t!]
\centering
\begin{tabular}{c}
    \includegraphics[trim = 0mm 0mm 0mm 0mm, clip, width=0.48\textwidth]{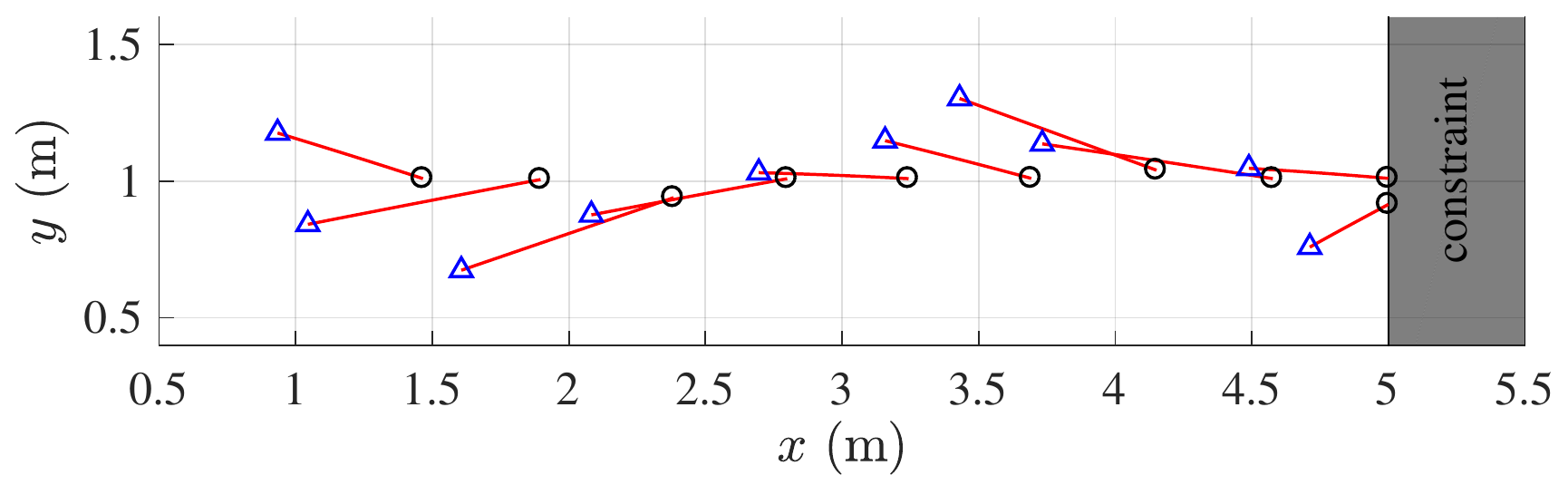}   \\
    \includegraphics[trim = 0mm 0mm 0mm 0mm, clip, width=0.45\textwidth]{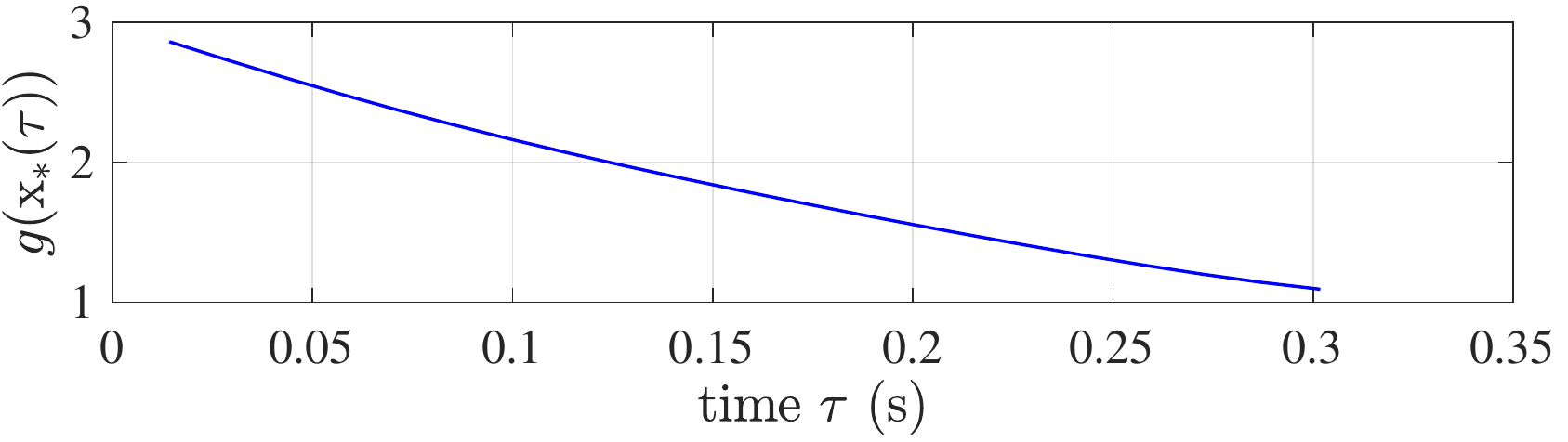} 
\end{tabular}
\caption{ (a) The red lines are the optimal state trajectories of the 10 robots. 
The blue triangles and the black circles are the position of the robots at the initial time and the optimal time, respectively.
Since the horizontal disturbance is big, some of robots on the right side cannot avoid the collision.
(b) The cost over $\tau$ ($g(\mathrm{x}_*(\tau))$) is minimized at $\tau=0.32$ (s). For $\tau>0.32$, the cost $g(\mathrm{x}_*(\tau))$ becomes the infinity since the state constraint is violated.
}
\label{fig:example2_traj}
\end{figure}

Figure \ref{fig:example2_traj} (a) shows the optimal state (position) trajectories of the robot system from the initial time to the minimum time ($\tau=0.32$ s) that minimizes the given cost.
Since the horizontal disturbance is more dominant than the robots' control, there exist some time ($\tau=0.32$ s) when one of the robots violates the state constraint in the given horizon. 
As shown in Figure \ref{fig:example2_traj} (b), the cost is minimized at $\tau=0.32$ (s) after which the state constraint is violated.

\subsection{Discussion}

Benefits of grid-based methods, such as the level set methods \cite{Osher02} and fast marching method \cite{sethian1996fast}, are that the global optimality is guaranteed, and a closed-loop control is provided. 
In other words, the optimal control is provided for any state-time pair.
However, due to computational complexity, grid-based methods are typically implemented offline to precompute the value function, and then used online with feedback.
Offline computation of grid-based methods is intractable for systems with continuous state dimension higher than six or seven.

The first and second examples in this paper are systems with forty and twenty state dimensions, whose dimensions are too high for grid-based methods.
Using the proposed method, the computation time for each example is 103.7 and 82.6 s, respectively.
In comparison to grid-based methods, this efficient computation is achieved since our method is gradient-based.
Although gradient-based methods provide local optimal solutions for non-convex problems, the proposed method guarantees the global optimality under the specified conditions in Table \ref{tab:comparisonCvxCond}.
Even though our methods provides an optimal control signal, unfortunately as with other gradient-based method, this is an open-loop control.

To have a closed-loop control, our method could cooperate with closed-loop control frameworks, such as model predictive control (MPC).
Closed-loop control is achieved by solving optimal control problems in real-time, and the real-time computation is allowed by the receding horizon setting.
Our method can cooperate with a number of MPC methods \cite{guanetti2018control} to have a closed-loop control and real-time computation.

\section{Conclusion and Future Work} 
\label{sec:conclusion}

This paper presents three Lax formulae for 1) Problem \ref{prob:1}, 2) Problem \ref{prob:2}, and 3) the time-invariant Problem \ref{prob:2}, and the Lax formulae are efficiently computed by gradient-based methods.
The derivation of the Lax formulae is based on proposed viscosity theory that provides a sufficient condition under which two different HJ equations have the same solution.
Furthermore, this paper provides a set of convexity conditions under which the Lax formulae provide an optimal solution by the gradient-based methods: for Problem \ref{prob:1}, the convexity in the state is required but not in the control; Problem \ref{prob:2} is generally non-convex; for the time-invariant Problem \ref{prob:2}, the dynamics and the stage cost only depend on the control, and the terminal cost and the state constraint are convex in the state.
This paper also presents a numerical algorithm for the Lax formulae.
For each class of problems, an example demonstrates how to utilize the Lax formulae and the numerical algorithm, and show the performance as well.

\appendix
\subsection{Proof of Theorem \ref{thm:ViscosityTheory}}
\label{appen:thm_ViscosityTheory}
\noindent

(i) The terminal values of $X_1$ and $X_2$ are the same.




(ii) $X_1$ is a subsolution to the PDE $0=F_2$. 

Since $X_1$ is a subsolution to the PDE $0=F_1(t,x,z,X_1,X_{1,t},D_x X_1,D_z X_1)$, for any $U\in C^1 ((0,T)\times\R^n\times\R)$ such that $X_1-U$ has a local maximum at $(t_0,x_0,z_0)\in(0,T)\times\R^n\times\R$ and $(X_1-U)(t_0,x_0,z_0)=0$,
\begin{align}
\begin{split}
    F_1(t_0,x_0,z_0,U_0,U_{t0},D_x U_{0},D_z U_{0})\geq 0,
    \label{eq:thm2_proof_1}
\end{split}
\end{align}
where $U_0=U(t_0,x_0,z_0)$, $U_{t0}=U_t(t_0,x_0,z_0)$, $D_x U_{0}=D_x U(t_0,x_0,z_0)$, and $D_z U_{0}=D_z U(t_0,x_0,z_0)$.
Since $X_1-U$ has a local maximum at $(t_0,x_0,z_0)$, $D_z U_0$ is in the superdifferential of $X_1$ with respect to $z$: $D_z U_0 \in D_z^+ X_1(t_0,x_0,z_0)$.
By \eqref{eq:ViscosityThm_Condition}, 
\begin{align}
    F_2(t_0,x_0,z_0,U_0,U_{t0},D_x U_{0},D_z U_{0})\geq 0.
    \label{eq:thm2_proof_2}
\end{align}

(iii) $X_1$ is a supersolution to the PDE $0=F_2$.

Since $X_1$ is a supersolution to the PDE $0=F_1$, for any $U\in C^1((0,T)\times\R^n\times\R)$ such that $X_1-U$ has a local minimum at $(t_0,x_0,z_0)\in[0,T]\times\R^n\times\R$ and $(X_1-U)(t_0,x_0,z_0)=0$,
\begin{align}
\begin{split}
    F_1(t_0,x_0,z_0,U_0,U_{t0},D_x U_{0},D_z U_{0})\leq 0,
    \label{eq:thm2_proof_3}
\end{split}
\end{align}
Since $X_1-U$ has a local minimum at $(t_0,x_0,z_0)$, $D_z U_0$ is in the subdifferential of $X_1$ with respect to $z$ ($ D_z^- X_1(t_0,x_0,z_0)$).
By \eqref{eq:ViscosityThm_Condition}, 
\begin{align}
    F_2(t_0,x_0,z_0,U_0,U_{t0},D_x U_{0},D_z U_{0})\leq 0.
    \label{eq:thm2_proof_4}
\end{align}
\qed

\subsection{Proof of Lemma \ref{lemma:cvx_V_z}}
\label{appen:cvx_V_z}
\noindent

(i) For all $y_1,y_2,y_3,y_4\in\R$,
\begin{align}
    \max\{ y_1 + y_2 , y_3 + y_4 \} \leq \max\{ y_1 , y_3 \} + \max\{ y_2 , y_4 \}.
    \label{eq:cvx_V_z_proof1}
\end{align}

(ii) Proof of \eqref{eq:cvx_V_z} for $V_1$.

Let 
\begin{align}
    &\text{Const}(\alpha,\tau) \coloneqq \max_{s\in[t,T]}c(s,\mathrm{x}(s)),\\
    &\text{Cost}(\alpha,\tau) \coloneqq \int_t^\tau L(s,\mathrm{x}(s),\alpha(s))ds + g(\tau,\mathrm{x}(\tau)).
\end{align}
\begin{align}
    & V_1(t,x, \theta_1 z_1 + \theta_2 z_2) \notag\\= &\min_{\alpha\in\mathcal{A}}\max_{\tau\in[t,T]}\max\big\{
    \theta_1 \text{Const}(\alpha,\tau) + \theta_2 \text{Const}(\alpha,\tau), \notag\\
    & \quad\quad\quad\quad\quad\textcolor{white}{i} \theta_1 \big[ \text{Cost}(\alpha,\tau) - z_1 \big] +\theta_2\big[ \text{Cost}(\alpha,\tau) - z_2 \big]
    \big\} \notag\\
    \leq &\min_{\alpha\in\mathcal{A}}\max_{\tau\in[t,T]} \theta_1 \max\big\{ \text{Const}(\alpha,\tau),  \text{Cost}(\alpha,\tau)-z_1\big\} \notag\\ 
    &\quad\quad\quad\quad + \theta_2 \max\big\{ \text{Const}(\alpha,\tau), \text{Cost}(\alpha,\tau)-z_2\big\}.
    \label{eq:cvx_V_z_proof2}
\end{align}
The first equality is according to the distributive property of the maximum operations, and the second inequality is by \eqref{eq:cvx_V_z_proof1}. 
For $\alpha\in\mathcal{A}$, we use $\tau_*(\alpha)$ to denote a maximizer of the last term in \eqref{eq:cvx_V_z_proof2}.
By the triangular inequality, we simplify \eqref{eq:cvx_V_z_proof2} to
\begin{align}
    &V_1(t,x,\theta_1z_1+\theta_2z_2)\notag\\
    \leq &\min_{\alpha\in\mathcal{A}} \theta_1 \max\big\{\text{Const}(\alpha,\tau),\text{Cost}(\alpha,\tau_*(\alpha))-z_1\big\} \notag\\ 
    +& \min_{\alpha\in\mathcal{A}} \theta_2 \max\big\{\text{Const}(\alpha,\tau), \text{Cost}(\alpha,\tau_*(\alpha))-z_2\big\} \notag\\
    \leq & \theta_1 V_1(t,x,z_1) + \theta_2 V_1 (t,x,z_2).
    \label{eq:cvx_V_z_proof3}
\end{align}
The last inequality holds by the definition of $V_1$ in \eqref{eq:ctrl_valueDef1}.

(iii) Similar to \eqref{eq:cvx_V_z_proof2},
\begin{align}
    & V_2(t,x, \theta_1 z_1 + \theta_2 z_2) \notag\\
    \leq &\min_{\alpha\in\mathcal{A}}\min_{\tau\in[t,T]} \theta_1 \max\big\{ \text{Const}(\alpha,\tau),  \text{Cost}(\alpha,\tau)-z_1\big\} \notag\\ 
    &\quad\quad\quad + \theta_2 \max\big\{ \text{Const}(\alpha,\tau), \text{Cost}(\alpha,\tau)-z_2\big\}.
    \label{eq:cvx_V_z_proof4}
\end{align}
Since the the last term in \eqref{eq:cvx_V_z_proof4} is greater than or equal to $\theta_1 V_2(t,x,z_1) +\theta_2 V_2 (t,x,z_2)$ by the triangular inequality, we conclude the proof.
\qed


\subsection{Proof of Lemma \ref{lemma:bound_gradient}}
\label{appen:lemma_bound_gradient}
\noindent

(i) The proof of \eqref{eq:lemma2_1}.

For $\bar{z} \geq 0$, $V_1(t,x,z+\bar{z}) \leq  V_1(t,x,z)$, and by the distributive property of the maximum operations,
\begin{align}
    V_1(t,x,z+\bar{z}) &=\inf_{\alpha\in\mathcal{A}}\max_{\tau\in[t,T]}\max\Big\{
    \max_{s\in[t,T]}c(s,\mathrm{x}(s))+\bar{z}, \notag\\
    &\int_t^\tau L(s,\mathrm{x}(s),u(s))ds+g(\tau,\mathrm{x}(\tau))-z
    \Big\}-\bar{z}\notag \\
    & \geq V_1(t,x,z) - \bar{z}. \notag
\end{align}
Thus, for $\bar{z}\geq 0$,
\begin{align}
    V_1(t,x,z) - \bar{z} \leq V_1(t,x,z+\bar{z}) \leq V_1(t,x,z),
    \label{eq:proof_subgrad_1}
\end{align}
and, by the same derivation, for $\bar{z}\leq 0$,
\begin{align}
    V_1(t,x,z) \leq V_1(t,x,z+\bar{z}) \leq V_1(t,x,z) - \bar{z}.
    \label{eq:proof_subgrad_2}
\end{align}

Suppose there exists $q>0$ in $ D_z^- V_1(t,x,z)$. Then, there exists $\epsilon>0$ such that $V_1(t,x,z+\bar{z})\geq V_1(t,x,z)+q\bar{z}$ for all $\bar{z}\in[-\epsilon,\epsilon]$. However, for $\bar{z}\in(0,\epsilon)$, 
\begin{align}
    V_1(t,x,z+\bar{z}) > V_1(t,x,z). 
\end{align}
This contradicts \eqref{eq:proof_subgrad_1}.

Suppose there exists $q<-1$ in $ D_z^- V_1(t,x,z)$. Then, there exists $\epsilon>0$ such that $V_1(t,x,z+\bar{z})\geq V_1(t,x,z)+q\bar{z}$ for all
$\bar{z}\in[-\epsilon,\epsilon]$. However, for $\bar{z}\in(-\epsilon,0)$, 
\begin{align}
    V_1(t,x,z+\bar{z}) > V_1(t,x,z) - \bar{z}. 
\end{align}
This contradicts \eqref{eq:proof_subgrad_2}.
Thus, $q\in[-1,0]$.

With the analogous derivation for $V_2$, we conclude that $ D_z^- V_2 (t,x,z)\subset[-1,0]$.

(ii) The proof of \eqref{eq:lemma2_2}.

The convexity of $V_i(t,x,z)$ ($i=1,2$) stated in Lemma \ref{lemma:cvx_V_z} implies that $D_z^+ V_i(t,x,z)$ contains a single superdifferential $q$ if $V_i(t,x,z)$ is locally affine in $z$, otherwise, $D_z^+ V_i(t,x,z)$ is the empty set. 
If $V_i(t,x,z)$ is locally affine in $z$, it is also differentiable in $z$.
As $\bar{z}$ converges to 0 in \eqref{eq:proof_subgrad_1}, we have $D_z V_i(t,x,z)\in[-1,0]$. Thus, if the superdifferential $q$ exists, $q$ and $D_z V_i(t,x,z)$ are the same, and $q\in[-1,0]$.
\qed


\subsection{Proof of Lemma \ref{lemma:Hamiltonian_comp}}
\label{appen:lemma_Hamiltonian_comp}
\noindent

(i) Case 1: $q=0$.
\begin{align}
\begin{split}
    \bar{H}(t,x,z,p,0) &= \max_{a\in A} -p\cdot f(t,x,a) = \max_{b\in B(t,x)} p\cdot b
\end{split}
\end{align}
by the definition of $B(t,x)$ in \eqref{eq:ctrlbound_B}.
Since $B(t,x)\subset \text{co}(B(t,x))$, 
\begin{align}
\begin{split}
    \bar{H}&(t,x,z,p,0) \leq \bar{H}_W(t,x,z,p,0).
\end{split}
\label{eq:lemma_Hamiltonian_comp_proof1}
\end{align}
On the other hand, let $b_*\in \arg\max_{b\in\text{co}(B(t,x))}p\cdot b$. Since $B(t,x)$ is compact, there exists a finite number of $b_i \in B(t,x)$ and $\theta_i\in[0,1]$ such that $b_* = \sum_i \theta_i b_i$ and $\sum_i \theta_i =1$. Then,
\begin{align}
     \bar{H}_W(t,x,z,p,0) &= \sum_i \theta_i p\cdot b_i 
     \leq \max_i \{p\cdot b_i\} 
     \notag \\
    & \leq \max_{b\in B(t,x)}p\cdot b = \bar{H}(t,x,z,p,0).
    \label{eq:lemma_Hamiltonian_comp_proof2}
\end{align}
The last inequality holds since all $b_i$s are in $B(t,x)$.
By \eqref{eq:lemma_Hamiltonian_comp_proof1} and \eqref{eq:lemma_Hamiltonian_comp_proof2}, we have
\begin{align}
    \bar{H}&(t,x,z,p,0)=\bar{H}_W(t,x,z,p,0).
\end{align}

(ii) Case 2: $q<0$.
\begin{align}
\begin{split}
    \bar{H}(t,x,z,p,q) &= \max_{a\in A} -p\cdot f(t,x,a)+qL(t,x,a) \\
    & = -q H\bigg(t,x, -\frac{p}{q} \bigg).
\end{split}
\label{eq:lemma_Hamiltonian_comp_proof3}
\end{align}
Since $H$ is convex in $p$ for each $(t,x)$ and lower semi-continuous in $p$, $H^{**}\equiv H$.
Thus, we have
\begin{align}
    \bar{H}_W(t&,x,z,p,q) 
     = -q \max_{b\in \text{co}(B(t,x))} -\frac{p}{q}\cdot b - H^*(t,x,b)\notag\\
    & = -q H^{**}\bigg(t,x,-\frac{p}{q}\bigg) = -q H\bigg(t,x,-\frac{p}{q}\bigg).
\label{eq:lemma_Hamiltonian_comp_proof4}
\end{align}
By \eqref{eq:lemma_Hamiltonian_comp_proof3} and \eqref{eq:lemma_Hamiltonian_comp_proof4}, we conclude 
\begin{align}
    \bar{H}&(t,x,z,p,q)=\bar{H}_W(t,x,z,p,q).
\end{align}
for all $q<0$.
\qed

\subsection{Proof of Lemma \ref{lemma:domain_ctrl_TI_Prob2}}
\label{appen:lemma_domain_ctrl_TI_Prob2}

This proof generalizes the proof of Lemma 1 \cite{lee2020hopf}, which is for the zero stage cost problem.

(i) For $b\in B(x)$, $b=-f(x,\bar{a})$ for some $\bar{a}\in A$. Then,
\begin{align}
    H_2^{\textnormal{TI}*}(x,b) &= \max_p -p\cdot f(x,\bar{a}) -H_2^{\textnormal{TI}}(x,p)\notag\\
    & \leq \max_p -p\cdot f(x,\bar{a}) +\min_{a\in A} p\cdot f(x,a) + L(x,a)\notag\\
    & < \infty.
    \label{eq:lemma_domain_ctrl_TI_Prob2_proof1}
\end{align}
The last inequality holds since $\min_{a\in A}p\cdot f(x,a)+L(x,a) \leq p\cdot f(x,\bar{a})+L(x,\bar{a})$ and $L$ is finite for a fixed $x$.

(ii) If $b=0$,
\begin{align}
    H_2^{\textnormal{TI}*}(x,b) =\max_p -H_2^{\textnormal{TI}} (x,p)\leq 0<\infty.
    \label{eq:lemma_domain_ctrl_TI_Prob2_proof2}
\end{align}

(iii) $b\in\text{co}(\{0\}\cup B(x))$

There exists a finite set of $\theta_i\in[0,1]$ ($\sum_i \theta_i \leq 1$), $a_i\in A$ such that 
    $b = -\sum_i \theta_i f(x,a_i).$
Since $H_2^{\textnormal{TI}*}$ is convex in $b$,
\begin{align}
    H_2^{\textnormal{TI}*}(x,b) \leq \sum_i \theta_i H_2^{\textnormal{TI}*}(x,b_i) + (1-\sum_i \theta_i)H_2^{\textnormal{TI}*}(x,0)<\infty
    \notag
\end{align}
by \eqref{eq:lemma_domain_ctrl_TI_Prob2_proof1} and \eqref{eq:lemma_domain_ctrl_TI_Prob2_proof2}.

(iv) $b\notin \text{co}(\{0\}\cup B(x))$

For the two convex sets $\{b\}$ and $\text{co}(\{0\}\cup B(x))$, by the separating hyperplane theorem \cite{boyd2004convex}, there exists a hyperplane ($P:\R^n \rightarrow \R$): $P (b') := p' \cdot b' + c$ such that $P(b)>0$ but $P(b')<0$ for all $b'\in\text{co}(\{0\}\cup B(x))$.
By picking $p= d p'$,
\begin{small}
\begin{align*}
    H_2^{\textnormal{TI}*} (x,b) \geq \max_{d} \min\big\{d p'\cdot b, \min_{a\in A} d p'\cdot(b+f(x,a))+L(x,a)\big\}.
\end{align*}
\end{small}
\negpar[-0em]Since $p'\cdot b>0$ and $p'\cdot(b+f(x,a))>0$ for all $a\in A$, the maximum of the right term in the above equation is attained at $d=\infty$, thus, $H_2^{\textnormal{TI}*} (x,b)=\infty$.
$\blacksquare$

\subsection{Proof of Lemma \ref{lemma:Hamiltonian_comp_TI}}
\label{appen:lemma_Hamiltonian_comp_TI}
\noindent
(i) Case 1: $q=0$.
\begin{align}
    \bar{H}^\textnormal{TI}_2 (x,z,p,0) = \max\{0,\max_{b\in B(x)}p\cdot b\},
\end{align}
where $B(x)$ is defined in \eqref{eq:ctrlBound_B_TI}. 
Since $B(x)\subset \text{co}(\{0\}\cup B(x))$,
\begin{align}
    \bar{H}^\textnormal{TI}_2 (x,z,p,0) \leq \bar{H}^\textnormal{TI}_W (x,z,p,0).
    \label{eq:lemma_Hamiltonian_comp_TI_proof1}
\end{align}
On the other hand, let $b_*\in\arg\max_{b\in\text{co}(\{0\}\cup B(x))}p\cdot b$, then, there exists a finite number of $b_i\in B(x)$ and $\theta_i\in[0,1]$ such that $b_*=\sum_i \theta_i b_i$ and $\sum_i\theta_i <1$. Thus, we have
\begin{align}
    \bar{H}_W^\textnormal{TI}(x,z,p,0) &= \sum_i \theta_i p\cdot b_i \leq \max\{0,\max_i p\cdot b_i\} \notag\\
    & \leq \bar{H}^\textnormal{TI}_2 (x,z,p,0).
    \label{eq:lemma_Hamiltonian_comp_TI_proof2}
\end{align}
Combining \eqref{eq:lemma_Hamiltonian_comp_TI_proof1} and \eqref{eq:lemma_Hamiltonian_comp_TI_proof2}, we have 
\begin{align}
    \bar{H}_W^\textnormal{TI}(x,z,p,0)=\bar{H}^\textnormal{TI}_2 (x,z,p,0).
    \label{eq:lemma_Hamiltonian_comp_TI_proof3}
\end{align}

(ii) Case 2: $q<0$.
\begin{align}
    \bar{H}_2^\textnormal{TI} (x,z,p,q)=-q \max\bigg\{0,H\bigg(x,-\frac{p}{q}\bigg)\bigg\}
    \label{eq:lemma_Hamiltonian_comp_TI_proof4}
\end{align}
and, since $H_2^{\textnormal{TI}**}\equiv H_2^{\textnormal{TI}}$ by the convexity of $H_2^{\textnormal{TI}}$ in $p$,
\begin{align}
    \bar{H}_W^\textnormal{TI} (x,z,p,q)=-q H_2^\textnormal{TI}\bigg(x,-\frac{p}{q}\bigg).
    \label{eq:lemma_Hamiltonian_comp_TI_proof5}
\end{align}
By combining \eqref{eq:lemma_Hamiltonian_comp_TI_proof4}, \eqref{eq:lemma_Hamiltonian_comp_TI_proof5}, and \eqref{eq:Ham_Ham2TI}, we conclude the proof. \qed
\subsection{Proof of Lemma \ref{lemma:cvx_vartheta1}}
\label{appen:lemma_conv_vartheta1_cont}
\noindent
\textbf{Proof.} \\
\noindent
(i) For each $k'\in\{0,...,K\}$, 
\begin{align}
    \sum_{k=0}^{k'} L(t_k, \mathrm{x}[k], \alpha[k])+g(t_{k'}, \mathrm{x}[k'])
\end{align}
is convex in $([\mathrm{x}[\cdot],\alpha[\cdot])$. 
Since a pointwise maximum of a family of convex functions is convex, the cost \eqref{eq:def_vartheta1_disc} is convex in $([\mathrm{x}[\cdot],\alpha[\cdot])$. 

(ii) The dynamical constraint, the control constraint, the initial state constraint, the state constraints are all convex in  $([\mathrm{x}[\cdot],\alpha[\cdot])$.
\qed
\subsection{Proof of Theorem \ref{thm:OptCtrl_PostProcess}}
\label{appen:thm_OptCtrl_PostProcess}
\noindent
\textbf{Proof.} Let 
\begin{align}
    & J_1(\tau) = \int_0^\tau H^*(s,\mathrm{x}_*(s),\beta_*(s))ds+g(\tau,\mathrm{x}_*(\tau)),\\
    & J_2(\tau) = \int_0^\tau L(s,\mathrm{x}_*^\epsilon(s),\alpha_*^\epsilon(s))ds + g(\tau,\mathrm{x}_*^\epsilon(\tau)).
\end{align}

Theorem 3 in \cite{donggunTAC_submitted} states that for any $\epsilon>0$, there exists $\delta>0$ such that, for any discretization $\{t_0=0,...,t_K=T\}$ where $\max_k |\Delta t_k|<\delta $:
\begin{align}
    &\| \mathrm{x}_* -\mathrm{x}_*^\epsilon \|_{L^\infty (0,T;\R^n)} < \epsilon, \label{eq:thm_OptCtrl_PostProcess_proof1}\\
    &\| J_1(\tau) - J_2(\tau) \|_{L^\infty (0,T)} < \epsilon.
    \label{eq:thm_OptCtrl_PostProcess_proof2}
\end{align}

(i) \eqref{eq:thm_OptCtrl_PostProcess_proof1} and Assumption \ref{assum:BigAssum} imply \eqref{eq:thm_OptCtrlProcess_eq1} and \eqref{eq:thm_OptCtrlProcess_eq2}. 

(ii) Let
\begin{align}
    \tau_*^1\in\arg\max_{\tau\in[0,T]} J_1(\tau), \quad \tau_*^2\in\arg\max_{\tau\in[0,T]} J_2(\tau).
\end{align}
Then, we have 
\begin{align}
    V_1(0,x) \geq J_1(\tau_*^2)~~\text{and}~~\max_{\tau\in[0,T]} J_2 (\tau) \geq J_2(\tau_*^1).
    \label{eq:thm_OptCtrl_PostProcess_proof3}
\end{align}
Since, for all $\epsilon>0$, there exists $\delta>0$ such that \eqref{eq:thm_OptCtrl_PostProcess_proof2} for both $\tau_*^1$ and $\tau_*^2$,
\begin{align}
    J_1(\tau_*^2) > \max_{\tau\in[0,T]} J_2 (\tau) -\epsilon, \quad J_2(\tau_*^1) > V_1(0,x) -\epsilon.
    \label{eq:thm_OptCtrl_PostProcess_proof4}
\end{align}
Combining with \eqref{eq:thm_OptCtrl_PostProcess_proof3} and \eqref{eq:thm_OptCtrl_PostProcess_proof4}, we conclude
\begin{align}
    |V_1(0,x)- \max_{\tau\in[0,T]} J_2 (\tau)| <\epsilon.
\end{align}
This concludes the proof of \eqref{eq:thm_OptCtrlProcess_eq3}. \qed

\subsection{Proof of Theorem \ref{thm:OptCtrl_PostProcess2}}
\label{appen:thm_OptCtrl_PostProcess2}
\noindent
\textbf{Proof.} Let 
\begin{align}
    & J_1(\tau) = \int_0^\tau H_2^{\textnormal{TI}*}(\mathrm{x}_*(s),\beta_*(s))ds+g(\mathrm{x}_*(\tau)),\\
    & J_2(\tau) = \int_0^\tau L(\mathrm{x}_*^\epsilon(s),\alpha_*^\epsilon(s))ds + g(\mathrm{x}_*^\epsilon(\tau)).
\end{align}

By applying Theorem 3 in \cite{donggunTAC_submitted}, Theorem 3 in \cite{lee2020hopf} is generalized: for any $\epsilon>0$, there exists $\delta>0$ such that, for any discretization $\{t_0=0,...,t_K=T\}$ where $\max_k |\Delta t_k|<\delta$:
\begin{align}
    & \| \mathrm{x}_*(\cdot) - \mathrm{x}_*^\epsilon(\sigma(\cdot;\beta_*^1))\|_{L^\infty (0,T;\R^n)} < \epsilon, \label{eq:thm_OptCtrl_PostProcess2_proof1}\\
    & \|J_1(\tau) - J_2(\sigma(\tau;\beta_*^1)) \|_{L^\infty (0,T)} <\epsilon. \label{eq:thm_OptCtrl_PostProcess2_proof2}
\end{align}
\eqref{eq:thm_OptCtrl_PostProcess2_proof1} and Assumption \ref{assum:BigAssum} imply \eqref{eq:thm_OptCtrlProcess2_eq1} and \eqref{eq:thm_OptCtrlProcess2_eq2}.
Also, \eqref{eq:thm_OptCtrl_PostProcess2_proof2} implies \eqref{eq:thm_OptCtrlProcess2_eq3}. \qed




\addtolength{\textheight}{0cm}   

\bibliographystyle{IEEEtran}
\bibliography{IEEEexample}

\vspace*{-2\baselineskip}
\begin{IEEEbiography}
    [{\includegraphics[width=1in,height=1.25in,clip,keepaspectratio]{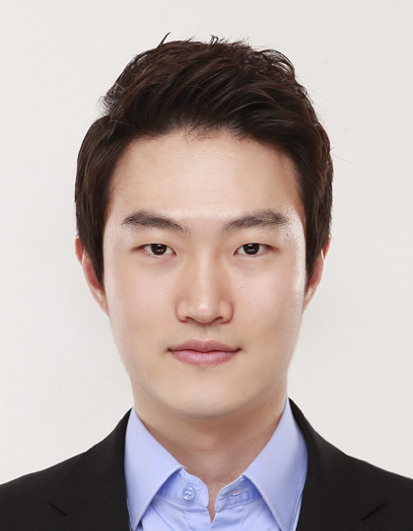}}]{Donggun Lee}
is a Ph.D. student in Mechanical Engineering at UC Berkeley. He received B.S. and M.S. degrees in Mechanical Engineering from Korea Advanced Institute of Science and Technology (KAIST), Daejeon, Korea, in 2009 and 2011, respectively. Donggun works in the area of control theory and robotics.
\end{IEEEbiography}
\vspace*{-2\baselineskip}
\begin{IEEEbiography}
    [{\includegraphics[width=1in,height=1.25in,clip,keepaspectratio]{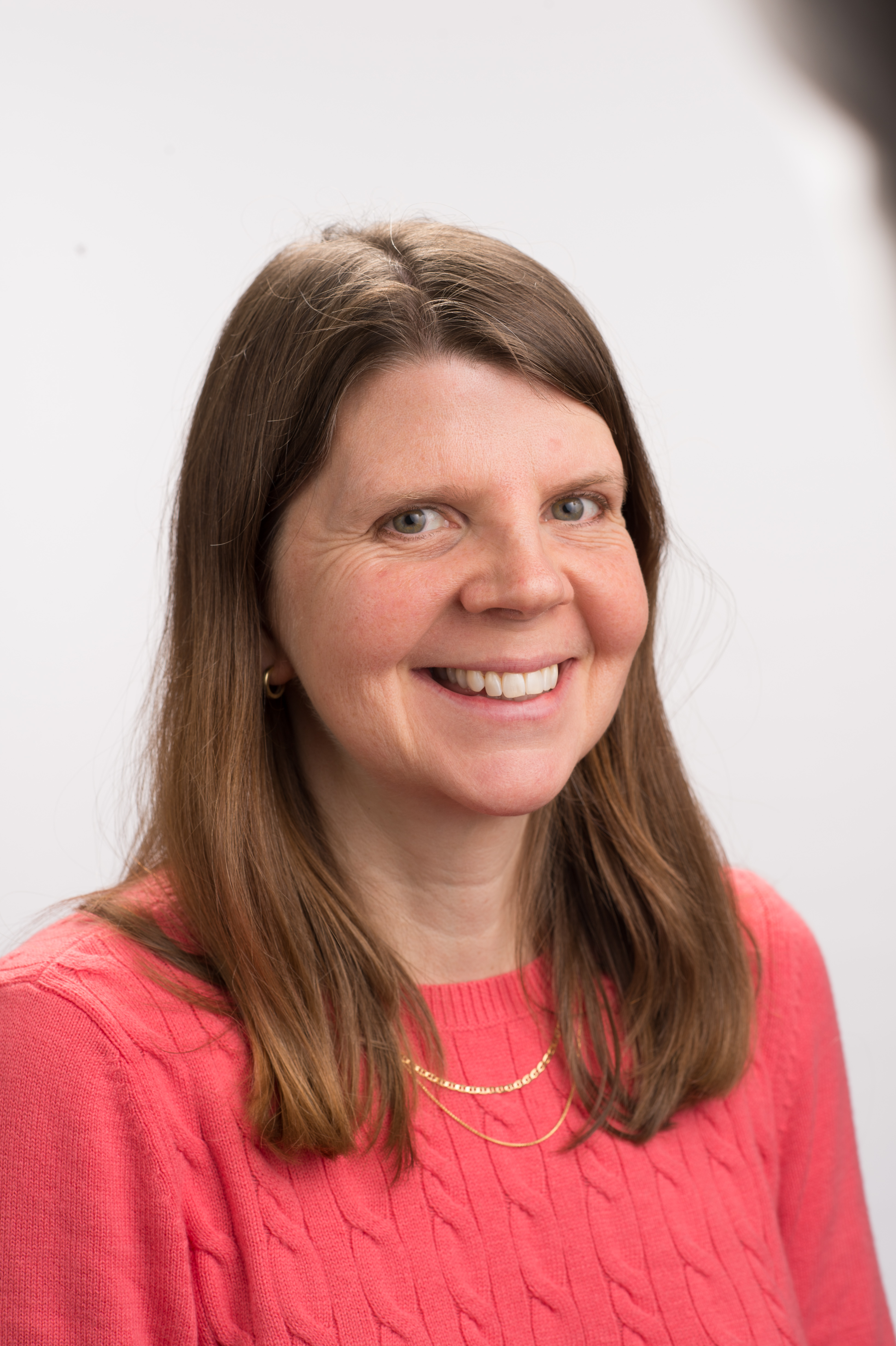}}]{Dr. Claire Tomlin}
is the Charles A. Desoer Professor of Engineering in EECS at Berkeley. She was an Assistant, Associate, and Full Professor in Aeronautics and Astronautics at Stanford from 1998 to 2007, and in 2005 joined Berkeley. Claire works in the area of control theory and hybrid systems, with applications to air traffic management, UAV systems, energy, robotics, and systems biology. She is a MacArthur Foundation Fellow (2006), an IEEE Fellow (2010), in 2017 she was awarded the IEEE Transportation Technologies Award, and in 2019 was elected to the National Academy of Engineering and the American Academy of Arts and Sciences.
\end{IEEEbiography}

\end{document}